\documentclass[12pt,
                headsepline=true,     
                footsepline=false]{article}

\usepackage[onehalfspacing]{setspace} 
\usepackage[a4paper, left = 2.5cm, right = 2.5cm, top = 2.5cm, bottom = 2.5cm]{geometry}

\usepackage{graphicx}%
\usepackage{multirow}%

\usepackage[title]{appendix}%
\usepackage{xcolor}%
\usepackage{textcomp}%
\usepackage{manyfoot}%
\usepackage{booktabs}%

\usepackage[utf8]{inputenc}
\usepackage[T1]{fontenc}
\usepackage[all]{nowidow}
\usepackage{graphicx}
\usepackage{csquotes}
\usepackage{amsmath, amsfonts, amssymb, bbm, bm}
\usepackage[normalem]{ulem}

\newtheorem{theorem}{Theorem}[section]
\newtheorem{cor}{Corollary}[section]
\newtheorem{proposition}{Proposition}[section]

\newtheorem{remark}{Remark}[section]%

\newtheorem{definition}{Definition}[section]%

\begin{document}

\title{Rectangular Gilbert Tessellation}
\date{}
\author{Emily Ewers\footnote{Mathematics Department, RPTU Kaiserslautern-Landau, Gottlieb-Daimler-Straße 48, 67663 Kaisers-lautern, Germany, emily.ewers@rptu.de}   \ and Tatyana Turova
\footnote{
Mathematical Center, University of
Lund, Box 118, Lund S-221 00, 
Sweden, \newline
tatyana@maths.lth.se}
}
\maketitle

\abstract{A random planar quadrangulation process is introduced as an approximation for certain cellular automata in terms of random growth of rays from a given set of points.
This model turns out to be a particular (rectangular) case of the well-known Gilbert tessellation, which originally models the growth of needle-shaped crystals from the initial random points with a Poisson distribution in a plane. 
From each point the lines grow on both sides of vertical and horizontal directions until they meet another line. 
This process results in a rectangular tessellation 
of the plane.
The central and still open question is the distribution of the length of line segments
in this tessellation.
We derive exponential bounds for the tail of this distribution. The correlations between the segments are 
proved to decay exponentially with the distance between their initial points. 
Furthermore, 
the sign of the correlation is investigated for some instructive examples.
In the case when the initial set of points is confined in a box $[0,N]^2$, it is proved that the average number of rays reaching the border of the box has a linear order in $N$.

}

\maketitle

\section{Introduction}
\setcounter{page}{1}

\subsection{From cellular automata to a planar quadrangulation}
The random quadrangulation process  considered here is motivated by a study \cite{ekstrom2022non}
of cellular automata introduced to describe neuronal activity.

Let us briefly recall the model from \cite{ekstrom2022non},  where a propagation 
of electrical impulses, or
activation, in a network on a  set of vertices $V\subseteq \mathbb{Z}^2$
 is described by a process $$(X_v(t), v\in V), \ t=0, 1, \ldots, $$ defined as follows.

The vertices in $V$ (which is typically a
finite set 
in a context of the neuronal models, but could be 
the entire $\mathbb{Z}^2$)
represent neurons, while the edges of the lattice $\mathbb{Z}^2$ are associated with synaptic connections. 
Hence,
each vertex $v$ is connected to its four nearest neighbours,
making an interacting neighbourhood of $v$:
\[N_v:=\Big\{
v + \{(0,-1), (0,1), (-1,0), (1,0)\} \Big\} \cap V.\]
 The set of vertices $V$ is split into two subsets:
 \begin{equation*}\begin{aligned}
 V^- &=\{(x,y)\in V: x,y \text{ both even}\}, \ \ \ \mbox{and}\\
 V^+ &= V\setminus V^-.
 \end{aligned}\end{equation*}
 All vertices in $V^-$ are assigned {\it inhibitory} type, while all 
vertices in $V^+$ are assigned {\it excitatory} type.

Let $X_v(t) \in\{0,1\}$ denote the state  of the process at vertex
$v{\in V}$ at time $t$:

$X_v(t) =1$ is the {\it active} state, and 

$X_v(t)=0$ is {\it nonactive}. 

\noindent
Given the set of initially active vertices
\[A(0) = \{v\in V : X_v(0)=1\}\subseteq V,\]
the dynamics of process $X(t)=(X_v(t), v \in V)$ in discrete time $t\geq 1$
is defined as follows:
\begin{equation}\label{T1*}
X{_v}(t+1)=\left\{
\begin{array}{ll}
1, & \mbox{ if } |N_v \cap V^+  \cap A(t)|-|N_v\cap V^-\cap A(t)|\geq 1, \\
0, & \mbox{ otherwise, }
\end{array}
\right.
\end{equation}
where
\begin{equation}\label{RA}
    A(t)=\{v\in V: X_v(t)=1\}
\end{equation}
is the set of active neurons (vertices) at time $t$.

In words, this definition reflects the properties of excitatory neurons to facilitate the propagation of  activity, while the inhibitory neurons, on the contrary, prevent the propagation of activity. These two actions result in non-monotonicity of the process $|A(t)|$, $t\geq 0$.
It also follows from definition (\ref{T1*}) that a neuron, which has more excitatory active neighbours than inhibitory active ones, becomes active, regardless of its own type.

The focus of study of this model is the evolution of the set of active vertices (\ref{RA}) for different initial sets $A(0)$. To underline the dependence on the initial set $A(0)=U$, we shall write
\begin{align}\notag
A(t)=A_U(t) = \{v\in V : X_v(t)=1 \mid A(0)=U\}.
\end{align}
By definition (\ref{T1*})
we have
for all $t\geq 0$
\begin{align}\notag
A_U(t+1) = 
\{v\in V : |N_v \cap V^+\cap A_U(t)|-|N_v\cap V^-\cap A{_U}(t)|\geq 1\}.
\end{align}
The question of interest concerns the limiting states (as $t\rightarrow \infty$) for such a process. Depending on the initial configuration $U$ of active vertices, 
qualitatively different 
limiting spatio-temporal states $A_U(t)$, called the activity patterns, may appear.  Their types are described in \cite{ekstrom2022non} as (i) sets expanding with a constant rate, (ii)
fixed contours of constant activity, or (iii) stable patterns moving across the sets of vertices $V$.
Notably, the lattice structure with its symmetries plays an important role here, particularly in the formation of the last two classes of states.

It must be admitted, however,  that a rigorous analysis of the model \cite{ekstrom2022non} is still beyond reach. Therefore, one may consider first some particular cases.
It is shown (\cite{ekstrom2022non}, Section 3.6) that for some sparse initial states $U$ (the set $U$ is sparse when for a vector $(x,y)$ between arbitrary two points in $U$, 
it holds that $\min\{|x|,|y|\}$ is larger than some  constant), the evolution of 
$A_U(t), \ 
t\in\{0,1, \ldots\},
$ 
can be described more straightforwardly as 
\[
A_U(t) = V \cap \widetilde{A}_U(t),
\]
where a random set 
$\widetilde{A}_U(t)\subset \mathbb{R}^2$ is monotone growing in continuous time $t\geq 0$ as follows.

\begin{definition}\label{dM}
Assume $N>2$,
\begin{equation}\label{ic}
    U \subset \{1,\ldots, N-1\}^2, 
\end{equation}
and
  let each  vertex $v\in U$ be assigned a random direction (vector $d_v$) independent of the rest: with equal probability $\frac{1}{2}$ this direction is either
 vertical ($d_v=(0,1)$) or horizontal ($d_v=(1,0)$). 
 Set $$\widetilde{A}_U(0):=U$$ to be the initial state at 
time $t=0$.
 At time $t=0$ from each
 vertex $v\in U$ two rays begin to grow along the edges of
the lattice $\mathbb{Z}^2$ into both opposite sides of the assigned directions ($+d_v$ and $-d_v$), 
with speed 1 (one edge per time unit). The growth of any ray stops
as soon as it meets another ray.
\end{definition}
\medskip

By this definition, in particular  for all $0<t<1$
the set $\widetilde{A}_U(t)$ is a union of closed intervals of equal length $2t$:
\begin{equation}\label{Rev1}
\widetilde{A}_U(t) = \cup_{v\in U} [v-t d_v, v+t d_v] ;
\end{equation}
here and below for any points $u,v \in \mathbb{R}^2$ we denote by $[u,v]$
the line segment between the end-points $u$ and $v$.

The assumption (\ref{ic}) on the initial activation set $U$ being confined in a bounded set $[1,N-1]^2$ implies that any ray which has crossed the boundary 
of set $V_N$ cannot meet on its way any other ray. 
Hence, after the time $N$ 
any ray will be either stopped by another ray within  $V_N$ or it will grow freely beyond $V_N$ towards infinity.
Therefore, the configuration of the rays within the bounded set $V_N$ will not change after (at most) time $N$, i.e. for all $t\geq N$
\begin{equation}\label{fg}
\widetilde{A}_U(t) \cap V_N  =\widetilde{A}_U(N) \cap V_N.
\end{equation}
The properties of the final set (\ref{fg}) are the object of study.

It is important to realize that Definition \ref{dM} reformulates the problem on a {\it non-monotone} evolution of cellular automata for a particular set of initial conditions
of the original model \cite{ekstrom2022non} in terms of a {\it monotone} growing process, which certainly simplifies the analysis. 

Note also that, due to the discrete structure of a lattice, the growth of rays in Definition \ref{dM} can be stopped by a collision of orthogonal rays at T-junctions or corners, but also by the \enquote{head-on} collision of rays with equal orientation. 
Whether all types of collisions can be realised in the course of the evolution of $\widetilde{A}_U$ depends, of course, on the initial set $U$.

To simplify this dependence, it is natural
to study,
as an approximation for the model described by
Definition \ref{dM},
a similar (and seemingly simpler) model of growing rays but on
$\mathbb{R}^2$. Then, assuming that the initial set $U$ is a random point process with a continuous distribution in $\mathbb{R}^2$, we $a.s.$ avoid the corner- and head-on
types of collisions of rays. 

It is remarkable that in reformulating our model in $\mathbb{R}^2$, we rediscover a well-known model by Gilbert, as we shall see in a moment.

\subsection{Model of planar quadrangulation}
We shall extend Definition \ref{dM} 
to describe a similar monotone growth 
of a set of rays in $\mathbb{R}^2$ 
from an arbitrary initial set of points $G(0)\subseteq \mathbb{R}^2$. For every $t>0$ the set $G(t)$ is a union of rays (including their endpoints) whose dynamics is defined as follows.

\begin{definition}\label{M1}
Let 
$\mathcal{V}$ be a locally finite set of points in $\mathbb{R}^2$ (i.e., within any finite square, there is a finite number of points of $\mathcal{V}$ ), and let a constant $0<p<1$ be fixed arbitrarily. For each $u\in \mathcal{V}$ 
a random vector $$d_u \in \{(0,1), (1,0)\}$$ is assigned independently of the rest with the probability
\begin{equation}\notag
\mathbb{P}\{d_u= (0,1) \}=p=1-\mathbb{P}\{d_u= (1,0) \}.
\end{equation}
Set 
$G(0):=\mathcal{V}.$
At time $t=0$, from every point $u\in \mathcal{V}$ into each side of the assigned direction ($\pm d_u$), a ray $\mathcal{L}^{\pm}_{u,d_u}(t)$
grows
with constant speed $1$ until it reaches another ray;
at this moment, the growth stops, and the ray remains unchanged thereafter. 

For any $t>0$, we define a collection of 
rays
\begin{equation}\label{Rev3}
G_{\mathcal{V}}(t):=  \left\{ \mathcal{L}^{+}_{u,d_u}(t),
\mathcal{L}^{-}_{u,d_u}(t),
u\in \mathcal{V}
\right\}.
\end{equation}
We call process $\left( G_{\mathcal{V}}(t), \ t\geq 0\right)$ 
the process of planar quadrangulation.
\end{definition}

We shall also consider the subset of 
$\mathbb{R}^2$, associated with the collection of rays $G(t)$, keeping (with a slight abuse) the above notation: 
\begin{equation}
G_{\mathcal{V}}(t) = 
\bigcup_{u\in \mathcal{V}} \left( \mathcal{L}^{+}_{u,d_u}(t)
\cup
\mathcal{L}^{-}_{u,d_u}(t)\right).
\end{equation}

   The case 
   when the initial set $ \mathcal{V}$ is confined to a bounded set is of particular interest here. 
   Let 
\begin{equation}\label{vN}
  \mathcal{V}_N   = \mathcal{V} \cap [0, N]^2.
\end{equation}
Then  the evolution of the set (\ref{Rev3}) within the same bounded square $[0, N]^2$, stops after a finite time. More precisely, at least for all $t>N$ it holds that 
   \begin{equation}\label{Evol}
       G_{\mathcal{V}_N}(t) \cap [0, N]^2 = G_{\mathcal{V}_N}(N) \cap [0, N]^2,
   \end{equation}
because by the time $t>N$ any ray with the initial point in $[0, N]^2$
either stopped growing or reached the boundary of the square $[0, N]^2$; in the latter case, it grows unbounded towards infinity. 
   
\begin{figure}[htbp]
    \centering
    \includegraphics[width = \linewidth]{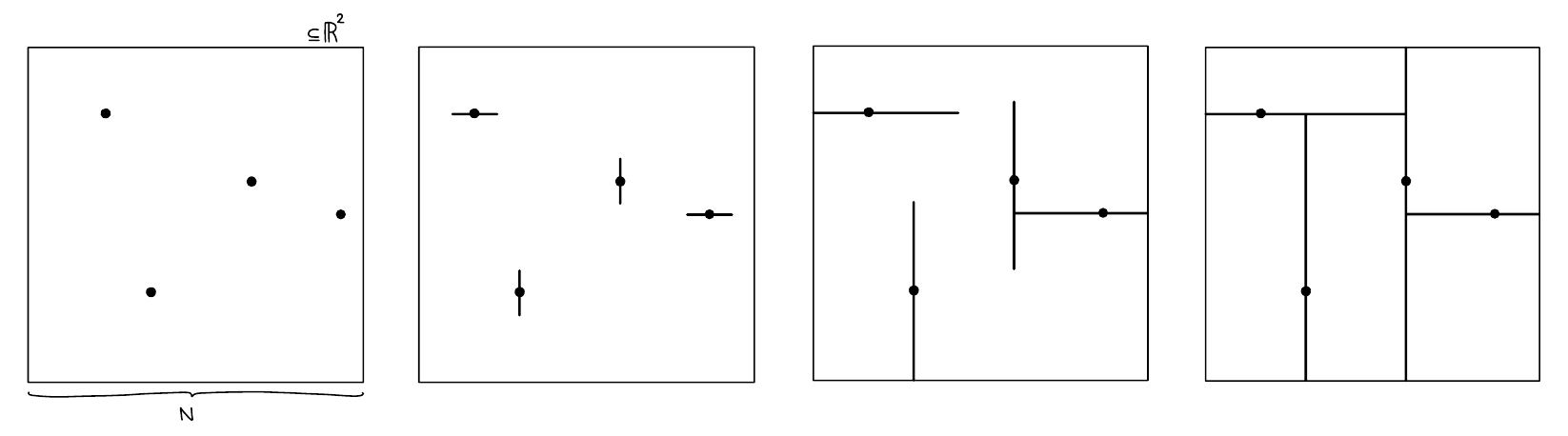}
    \caption{Example of a process with four points from the left to the right: (1) at $t=0$, (2) at some moment before any collision, (3)
    at the time of the first collision, (4) at $t=N$.}
    \label{fig:setup}
\end{figure}

   The process is illustrated in Figure \ref{fig:setup}.

\begin{figure}[htbp]
    \centering
    \includegraphics[width=0.6\linewidth]{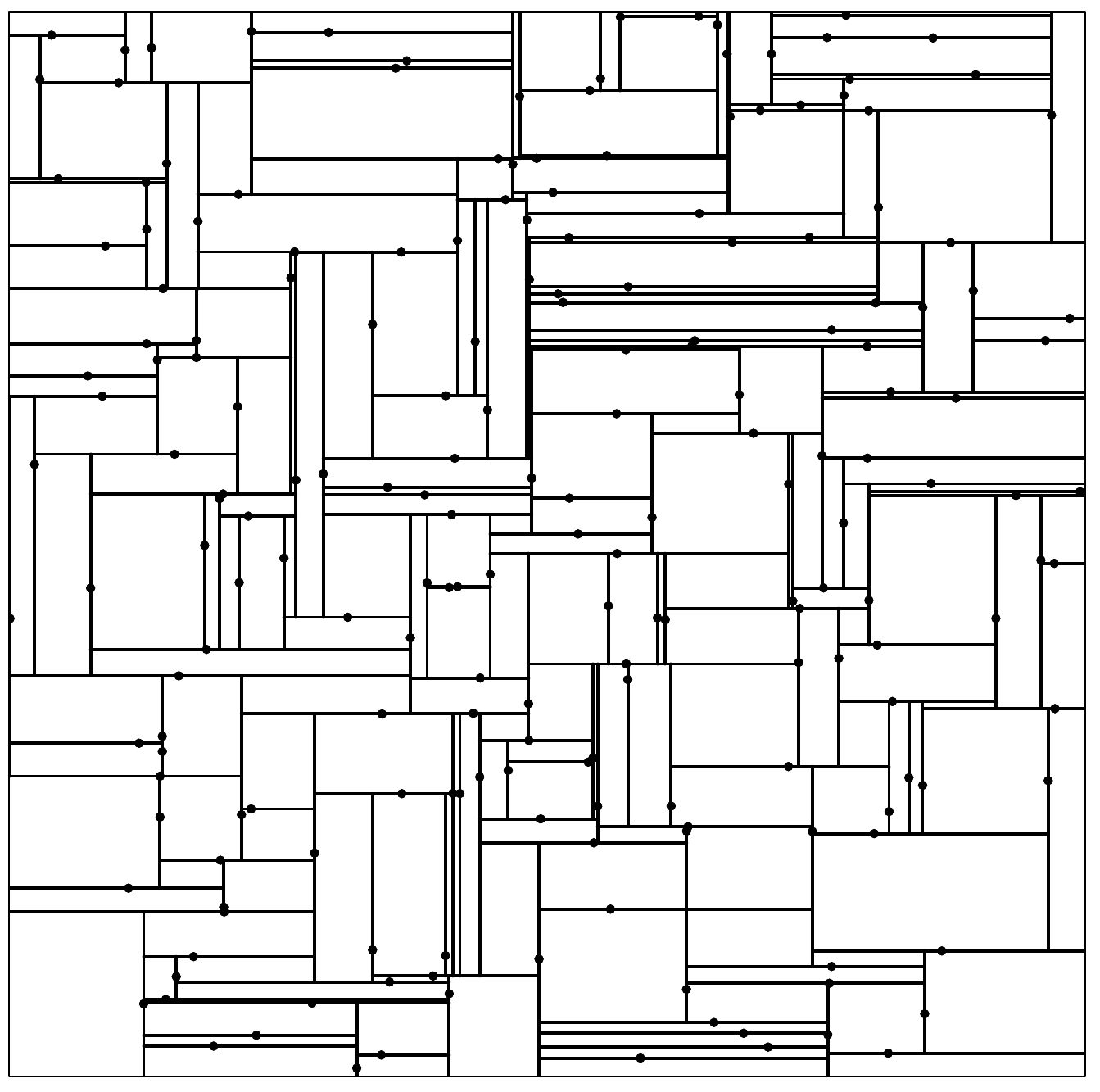}
    \caption{Example of a quadrangulation with 200 initial points.}
    \label{fig:intro}
\end{figure}

When all the rays within a square $[0, N]^2$ stop growing,
   all faces in the resulting configuration are almost surely rectangles, and moreover, all the junctions of the rays are $T$-shaped.
   An example of such a quadrangulation of a square is shown in Figure \ref{fig:intro}.

In the sequel, we assume that the initial set $\mathcal{V}$ is random; it is a realization of a 
   (marked) stationary Poisson point process $\eta$ in $\mathbb{R}^2$ with constant intensity \mbox{$\lambda>0$}. Each point 
   $u$
   is equipped with a random independent direction $d_u$. The latter means that 
   given a non-empty realization $U$ of the random set 
    $\mathcal{V}$,
   for each point  $u\in U$, a random direction $d_u$ is assigned independently of other points' directions, and of the location of the set $U$.
   By this assumption, for example, the number of points in $[0,N]^2$  follows a Poi$(\lambda N^2)$ distribution.
   
   The initial conditions, i.e., the set $\mathcal{V}$ and the directions $d_u$, $u\in \mathcal{V}$,  are the only source of randomness; otherwise, the dynamics are deterministic.
   In other words, given a realization  of the set $\mathcal{V}$
   equipped with direction vectors,   the process is deterministic.

\subsection{Gilbert tessellation}\label{GT}

It is remarkable that the model in Definition \ref{M1} turns out to be a particular case of the Gilbert tessellation, introduced already in 1967  by E.N. Gilbert \cite{Gilbert1}. In fact, originally, Gilbert defined and studied a much more general model on $\mathbb{R}^D$ in \cite{Gilbert2}, inspired by the facts on crystals growing in a metal or mineral (if $D=3$). Furthermore, Gilbert wrote \enquote{The
\enquote{crystals} may also be cells in living tissue}
\cite{Gilbert2}, to which statement the model \cite{ekstrom2022non}
described in the introduction
fits particularly well.

In \cite{Gilbert1} the Gilbert tessellation model is considered on $\mathbb{R}^2$, and
it is assumed that the directions ${d}_u$ are independent and uniformly distributed on $[0, 2\pi]$. 
Specifically, the rectangular case of Gilbert tessellation was first investigated in \cite{MackisackMiles1996FullGilbert} where Gilbert's mean-field approach to study the distribution of the length of rays was further developed.
 
 The approximation results of \cite{MackisackMiles1996FullGilbert} 
 for the rectangular Gilbert tessellation
 were further improved in
 \cite{BurridgeCowan2016HalfGilbert} and \cite{BurridgeCowanMa2013HalfGilbert};
a certain modification was introduced, the so-called \enquote{Half Gilbert tessellation} model for which a new technique was developed to get exact results on the distribution of the length of rays. It is argued in \cite{BurridgeCowanMa2013HalfGilbert}
 that the results of the modified model give a good approximation for the original model.
Still, no exact results are available on the expected value of the line segments' length in the original Gilbert tessellation. 

We take a different approach to study the original rectangular model. Namely, we derive bounds for the tail of the segment length distribution. Obviously, this does not yield the exact expectation, but the result allows us to study the correlations in the model. 
It can be noted that the derived exponential bounds 
for the tail are in line with the related limit theorems 
and the 
 exponential fast stabilization of the Gilbert tessellation
 proved in \cite{S}. 
 
 Very recently a new
generalisation of Gilbert tessellation, namely \enquote{Iterated Gilbert mosaics}, was introduced in \cite{BT}. The latter work opens a new
perspective of Gilbert tessellation 
from a point of view of tropical geometry.

Our results on the quadrangulation process are stated in Section 2. The proofs are provided in Section 3.

\section{Results}

\subsection{Rectangular Gilbert tessellation model}

From now on we concentrate on the
model $G_{\mathcal{V}}(t)$ given by Definition \ref{M1} with $p=\frac{1}{2}$, 
and random set $\mathcal{V}$ generated by a Poisson process $\eta$ in $\mathbb{R}^2$ with intensity $\lambda$;
we call it {\bf rectangular Gilbert tessellation model}.

By the definition the initial state  is  a random set of points $G_{\mathcal{V}}(0)=\mathcal{V}$. Given a set $\mathcal{V}$,  each point $u \in \mathcal{V}$ is equipped with a random independent direction $d_u$ (or a \enquote{mark}, would be a typical name in the context of point processes, see e.g. \cite{Kingman}), horizontal or vertical, equally probable (as $p=\frac{1}{2}$ in Definition \ref{M1}).
Given  set $\mathcal{V}$ of points with assigned directions $d_u, u \in \mathcal{V}$, at any time $t\geq 0$ the state of the process $G_{\mathcal{V}}(t)$ by  the Definition \ref{M1} is a set of closed line segments (i.e., the rays with their endpoints):
\begin{equation}\label{Rev4}
G_{\mathcal{V}}(t) = 
\bigcup_{u\in \mathcal{V}} \left( \mathcal{L}^{+}_{u,d_u}(t)
\cup
\mathcal{L}^{-}_{u,d_u}(t)\right),
\end{equation}
where for all $u\in \mathcal{V}$
\[\mathcal{L}^{\pm}_{u,d_u}(0)=u.\]
Recall again that given the initial set of points equipped with directions, the dynamics of the process is deterministic.

\subsection{Bounds for distribution of the rays' length}
The main question of interest is the limiting distribution (when $T\rightarrow \infty$) of the length of line segments $$|\mathcal{L}^{+}_{v,d_v}(T)
\cup
\mathcal{L}^{-}_{v,d_v}(T)|=|\mathcal{L}^{+}_{v,d_v}(T)
|+|
\mathcal{L}^{-}_{v,d_v}(T)|
, \ v\in \mathcal{V}.
$$

Recall that the characteristic
    property of the Poisson process is independence for the processes over non-intersecting areas (see \cite{Kingman}, \cite{Cox}).
    Therefore, for any $v\in \mathbb{R}^2$ fixed arbitrarily, conditionally on $v\in \mathcal{V}$, i.e., that point $v$ belongs to the Poisson process $\eta$,
    the distribution of the points of $\eta$ in $\mathbb{R}^2\setminus \{v\}$ is still Poisson with the same intensity $\lambda$. 
Therefore, we may fix $v\in \mathbb{R}^2$ arbitrarily and given also $d_v$, consider the growth of the rays $\mathcal{L}^{\pm}_{v,d_v}$ from $v$ in the random environment of the Poisson process $\eta$.

First we observe that by the Definition \ref{M1}
\begin{equation} 
  {L}_{v,d_v}^\pm (T):=|\mathcal{L}^{\pm}_{v,d_v}(T)|\leq T 
\end{equation}
for all $T>0$, as the speed of growth of the rays is one.
(Here and below we use the upper index $\pm$ if the statement equally concerns both cases, $+$ and $-$, separately.)
For the same reason for any $t<T$ the event 
\begin{equation} \label{lt2}
  \{{L}_{v,d_v}^\pm (T)\geq t \}
\end{equation}
implies that the ray 
$\mathcal{L}^{\pm}_{v,d_v}$ still grows at time $t$, which is the event
\begin{equation} \label{lt1}
   \{{L}_{v,d_v}^\pm (t)=t \};
\end{equation}
otherwise, it could not reach the length $t$, since once the growth stops, it stops forever.
On the other hand, event (\ref{lt1})
yields the event (\ref{lt2}). 
Hence, for any $T>t>0$ we have 
\begin{equation} \label{lt}
  \{{L}_{v,d_v}^\pm (T)\geq t \}=   \{{L}_{v,d_v}^\pm (t)=t \}
\end{equation}
(consider Figure \ref{fig:setup} to gain some intuition).

Notice that the distance and time could be interchangeable in our model due to a constant speed of growth (equal to $1$). Therefore, we mostly use the variable $t$ for the distance as well.

\begin{theorem}
\label{T1} Let  $G_{\mathcal{V}}(t)$ 
be a rectangular Gilbert tessellation model with 
random set $\mathcal{V}$ generated by a Poisson process in $\mathbb{R}^2$ with intensity $\lambda$. 

There are positive constants 
    $c_1, c_2$ and $\alpha_1 \geq \alpha_2$ such that for 
    any $\lambda >0$, for any
    arbitrarily fixed $v\in \mathbb{R}^2$ and $d_v\in \{(1,0), (0,1)\}$
    \begin{equation}\notag
      c_1 e^{-\alpha_1 {\sqrt{\lambda}}t} \leq  \mathbb{P}\{
      |\mathcal{L}^{+}_{v,d_v}(T)
\cup
\mathcal{L}^{-}_{v,d_v}(T)|
      >t \mid v\in\mathcal{V}\}\leq c_2 e^{-\alpha_2 {\sqrt{\lambda}}t}
    \end{equation} 
    for any fixed $t>0$ uniformly in  $T>t$.
\end{theorem}

This result tells us that 
the total length of the line segment through point $v$ in  $G_{\mathcal{V}}(T)$ has an exponentially decaying tail of distribution when $T\rightarrow \infty$.
Whether the limiting distribution is indeed exponential is not resolved by this result, although it shows a correct scaling for the expectation which is $1/\sqrt{\lambda}$. The latter simply reflects the following scaling property of the model. To underline the dependence on $\lambda $, let us write $L^{\pm}_{v, d_v}(T)= L^{\pm}_{v,d_v, \lambda}(T)$. Then it follows directly by the properties of the Poisson process 
 that 
 \begin{equation*}
    L_{v, d_v,\lambda}^{\pm}(T)\stackrel{d}{=}\frac{1}{\sqrt{\lambda}}L_{v, d_v,1}^{\pm}(\sqrt{\lambda}T).
\end{equation*}

\subsection{Correlations}
The core of difficulties in the analysis of Gilbert tessellation is in the spatial dependence. Therefore
 we begin to investigate the covariance structure of the model, which has not been done previously. 
 
First, we characterise some independent events in the model with the help of the following construction.

\begin{figure}[htbp]
    \centering
    \includegraphics[width = 0.5\linewidth]{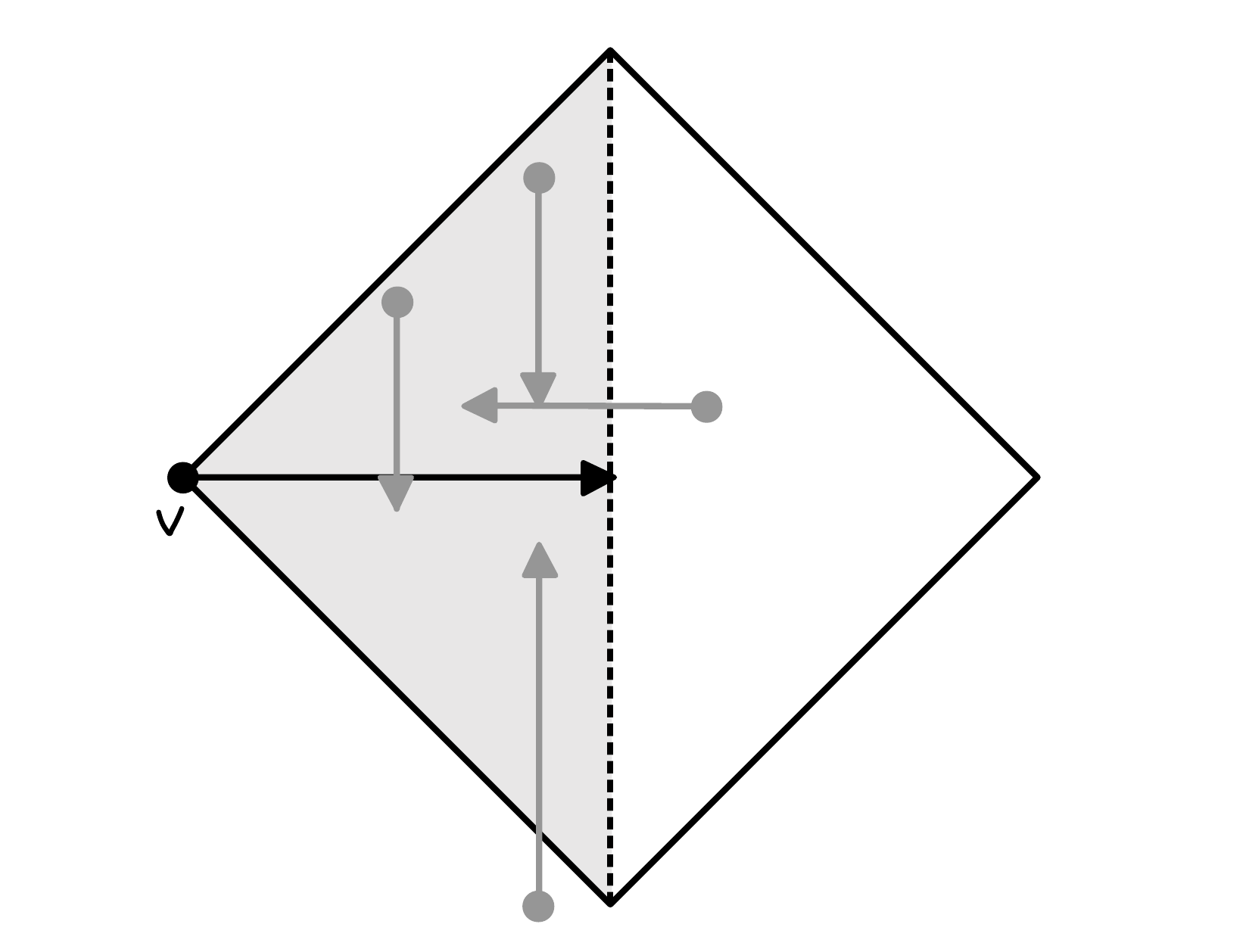}
    \caption{The shaded area is the \textit{cone} $C_{v,d_v}^+(t)$, 
    the entire square is  $D_{v,d_v}^+(t)$, including examples of other lines growing and influencing a potential ray ${\cal L}_{v,d_v}^+(t)$ marked by the right arrow pointed to the centre of the square.}
    \label{fig:coneDiamond}
\end{figure}

  Given $v\in\mathcal{V}$ and a direction $d_v$ denote
  $D_{v,d_v}^+ (t)$ (and $D_{v,d_v}^- (t)$) for $t>0$ a square with a corner in $v$, whose  diagonal from this corner 
  is the interval $[v,v+2td_v]$ (correspondingly, $[v-2td_v, v]$).
  Let $C_{v,d_v}^\pm (t)$ denote the half of this square, which is the (2-dimensional) cone with its  apex at $v$; see Figure \ref{fig:coneDiamond}.
Hence, a ray $\mathcal{L}_{v,d_v}^\pm(t)$
would follow the centre line of $C_{v,d_v}^\pm (t)$.

\noindent
\begin{definition}\label{depA}
 For all $t>0$ 
 we call 
 $D_{v,d_v}^+ (t)$ the area of dependence for the event 
 $$A_{v,d_v}^+(t):=\{
 {L}_{v,d_v}^+(t) =t\}.$$
 Correspondingly, for all $t>0$ 
 we call $D_{v,d_v}^- (t)$ the area of dependence for the event $$A_{v,d_v}^-(t):=\{{L}_{v,d_v}^-(t) =t\}.$$
\end{definition}

By the definition of our process,
a ray, which might stop 
the growth of ray 
$\mathcal{L}_{v,d_v}^\pm$ within time $t$, must have its initial point in the area
$C_{v,d_v}^\pm (t)$, while each of such rays, in turn, can be stopped within time $t$
by another ray only if the latter is 
originated in the area 
$D^\pm_{v,d_v} (t)$. This observation leads to the following useful statement on independence, checked by the same straightforward argument.

\begin{proposition}\label{Prop1}
  For all $t>0$ and $ u\neq v$, 
the events $A_{u,d_u}^\pm(t)$ and $A_{v,d_v}^\pm(t)$,  conditionally on $\{u,v \in \mathcal{V}\}$, are
  independent if the corresponding
 areas of dependence do not intersect, i.e., if
 \begin{equation}\label{CI}
 D_{u,d_u}^\pm (t) \cap D_{v,d_v}^\pm (t)=\emptyset.
 \end{equation}
\end{proposition}

Proposition \ref{Prop1} immediately yields the independence of the rays growing from one point but  in different directions.

\begin{cor}\label{CorI}
    Given 
    $v\in\mathcal{V}$ 
    and the direction $d_v$ the variables 
    $L^+_{v,d_v} (t)$ and $L^-_{v,d_v}(t)$ are independent for all $t$.
    
    \end{cor}

Let us use here the following definition of distance.
\begin{definition}\label{def:distVert}
    For any  $u=(u_1,u_2),v=(v_1,v_2)\in \mathbb{R}^2$ set
    \begin{equation*}
        \Vert u-v\Vert_1=|u_1-v_1|+|u_2-v_2|.
    \end{equation*} 
\end{definition}
Observe that for two different points  $u,v$ given that  $u,v \in \mathcal{V}, $
the areas  $D_{u,d_u}^\pm (t)$ and $  D_{v,d_v}^\pm (t)$ of dependence for the corresponding rays   do not intersect at least up to time $t=\frac{1}{4}\Vert u-v\Vert_1.$ Note that the latter is the minimal time $\frac{1}{4}\Vert u-v\Vert_1$ until
$D_{u,d_u}^\pm (t)$ and $  D_{v,d_v}^\pm (t)$  meet (intersect), and it is
indeed achieved for example, when $d_u\neq d_v$ and $$3|u_i-v_i|=|u_j=v_j|, \ \ i\neq j.$$
Hence, 
Proposition \ref{Prop1} implies in this case the following.

\begin{cor}\label{le:distForInd}
Let two points $u,v  \in \mathbb{R}^2$ be fixed arbitrarily. Given that \mbox{$\{u,v \in \mathcal{V}\}$}, and given 
any direction vectors $d_u$ and $d_v$
the events 
$A_{u,d_u}^\pm(t)$ and $A_{v,d_v}^\pm(t)$ 
 are conditionally independent for all $ t\leq \frac{1}{4}\|u-v\|_1$.
\end{cor}

As a measure of dependence between
two events $A$ and $B$ (conditionally, perhaps  on another event $C$) we shall consider the covariance between their corresponding indicators:
\begin{equation}\label{CABu}
     {\bf C}\left( A,B \right)
:=
\mathbb{P}\{A\cap B\} 
-\mathbb{P}\{A\} \mathbb{P}\{ B\},
\end{equation}
or
\begin{equation}\label{CAB}
     {\bf C}\left( A,B \mid C\right)
:=
\mathbb{P}\{A\cap B\mid C\} 
-\mathbb{P}\{A\mid C\} \mathbb{P}\{ B\mid C\}.
\end{equation}

Generally speaking, the correlation between the  events 
$A_{u,d_u}^\pm(t)$ and $A_{v,d_v}^\pm(t)$ is involved,  and even its sign does depend on the distances between the origins of the rays, on the overlap of the areas of dependence, and on the directions of the rays. 
Below, we shall provide some examples.

\subsubsection{Example of negative correlations}
In some cases, the correlation sign is straightforward to establish.
In particular, it is easy to find a pair of points $u\neq v$ together with the direction vectors $d_u, d_v$, such that
\begin{equation}\label{S2}
\mathbb{P}\{
A_{u,d_u}^+(t) \cap A_{v,d_v}^+(t)
\mid u,v  \in\mathcal{V}\}=0
\end{equation}
for  some $t\geq 0$.
Indeed, 
choose $v=(v_1,v_2)$ 
and $u=(u_1,u_2)$
with $v_1<u_1$
and $u_2<v_2$, and let $d_v=(1,0), d_u=(0,1)$. Then for any 
$$t>\|u-v\|_1$$
we have (\ref{S2}), 
which yields in this case
\[
{\bf C}\left(A_{u,d_u}^+(t), A_{v,d_v}^+(t)\mid u,v  \in\mathcal{V} \right)<0.
    \]

 On the other hand, if two points in $\cal{V}$ are close enough, then the  rays growing from these points in the same direction are positively correlated, as we see below.

\subsubsection{Example of positive correlations}
\begin{proposition}\label{Prop2}
    Let $u,v\in \mathbb{R}^2$ be an arbitrary pair of points with  equal assigned directions \mbox{$d_u=d_v$}. 
    
    There is a constant 
    $c_3$ such that if 
    \[\|u-v\|_1<c_3\frac{e^{-\alpha_1 {\sqrt{\lambda}}t}}{\lambda t}\]
    ($\alpha_1$ is from Theorem \ref{T1}),
    then 
    \begin{equation}\notag
    {\bf C}\left(A_{u,d_u}^+(t), A_{v,d_v}^+(t)\mid u,v  \in\mathcal{V} \right)
    >0
    \end{equation}
    for all large 
    $t$. (A similar statement holds when the \enquote{$+$}-direction is replaced by the \enquote{$-$}-direction.)
\end{proposition}
The proof of this proposition requires results from the previous section, and it will be provided below.
\subsubsection{Exponential decay of correlations}
Finally, we shall derive the exponential bound for the decay of the absolute value of correlations (hence, no matter the sign) between the events 
$A_{u,d_u}^{\pm}(t)$
and $A_{v,d_v}^{\pm}(t)$ with respect to the distance between $u$ and $v$.

\begin{theorem}\label{th:covExpBound}
Let $v,u\in \mathbb{R}^2$ be an arbitrary pair of points with  given arbitrary directions \mbox{$d_u, d_v$}. 
 There are positive constants $c$ and $\alpha$ such that for any $\lambda$
\begin{equation}\label{eq:covDec}
    \left|{\bf C}\left(A_{u,d_u}^{\pm}(t), A_{v,d_v}^{\pm}(t)\mid u,v  \in\mathcal{V} \right)\right|
\leq c e^{-{\alpha}{\sqrt{\lambda}}\|u-v\|_1}
    \end{equation}
    for all $t>0$.
\end{theorem}
Observe, that in view of Corollary \ref{le:distForInd} the exponential decay established in the last theorem is meaningful when $ t> \frac{1}{4} \|u-v\|_1$; otherwise, the statement (\ref{eq:covDec}) is trivial.

\subsection{Bounded initial state}
 Here we consider a tessellation  model $G_{\mathcal{V}_N}(t)$ where
    the initial set $\mathcal{V}_N$  is confined to a bounded set $\Lambda_N:=[0,N]^2$
    as defined in (\ref{vN}):
   \[\mathcal{V}_N   = \mathcal{V} \cap \Lambda_N.\] 
As it was established in (\ref{Evol}), the configuration of rays in $G_{\mathcal{V}_N}(t)$
within $\Lambda_N$ does not change for $t>N$.

\subsubsection{Quadrangulation of a square}

   First, we shall derive a
rather curious 
 exact result
about the number of rectangles in the square $\Lambda _N$ induced by this quadrangulation; recall (\ref{Evol}).

\begin{proposition}\label{Plg} Let $N>2$ be arbitrarily fixed. 
     Given the value $|\mathcal{V}_N|=|\mathcal{V} \cap \Lambda_N|$, 
    the number of rectangles in the planar graph constructed 
  by the line segments of the set 
  \begin{equation}\label{Fcon}
  G_{\mathcal{V}_N}(N)\cap \Lambda_N 
  \end{equation}
  together with the boundary of the square $\Lambda_N$,
    equals $a.s.$
    \begin{equation}\label{Euler}
        |\mathcal{V}_N|+1.
    \end{equation}
\end{proposition}
\noindent
{\it Proof.} For each realization $\mathcal{V}_N$,
the number of vertices of the planar graph corresponding to the final configuration (\ref{Fcon}) is 
\[V=2|\mathcal{V}_N | +4,\]
where $4$ is the number of corner vertices of the square $\Lambda_N$, and $2|\mathcal{V}_N |$ vertices are at $T$-junctions. Notice that here only the points of $T$-junctions and the corner points of set $[0,N]^2$ itself are the proper {\it vertices of the planar graph}; not to be confused with the points of the Poisson process.

Counting all corners of the rectangular faces created by these vertices and dividing the result by $4$ gives us the number of rectangles, which is
$$\frac{(2|\mathcal{V}_N |)2+4}{4}=|\mathcal{V}_N |+1,$$
as stated in (\ref{Euler}).\hfill$\Box$

We can also derive the number of edges of the planar graph
corresponding to the final configuration (\ref{Fcon}). Note that the edges here are the segments between the vertices of the planar graph, and not simply the sides of the rectangles. 
\begin{cor}
    The number of edges of the planar graph defined in Proposition \ref{Plg}
corresponding to the final configuration equals $a.s.$
    \begin{equation}\label{Euler*}
        3|\mathcal{V}_N|+4.
    \end{equation}
\end{cor}
\noindent
{\it Proof.}
The Euler’s formula, which binds the number
of vertices ($V$), edges ($E$) and faces ($F$) for planar graphs, states: $$V-E+F=2.$$ 
Adding to the number of rectangles one outer face gives us the total number of faces:
\[F=|\mathcal{V}_N |+2.\]
Now by the Euler's formula
$$V-E+F= (2|\mathcal{V}_N | +4)-E+(|\mathcal{V}_N |+2)=2.$$
This yields the result (\ref{Euler*}). \hfill$\Box$

\subsubsection{Percolation}

Observe that the final quadrangulation of $\Lambda_N$ 
by $G_{\mathcal{V}_N}(N)$
as discussed above does not change when the growth of lines continues beyond the boundary of the square $\Lambda_N$. 
Nothing prevents the growth of the rays outside of $\Lambda_N$, so all the rays which reach the boundary will grow unbounded towards infinity, 
i.e., \enquote{percolate}, correspondingly in four possible directions. On the other hand, these lines do not add anything inside the square $\Lambda_N$.

A phenomenon of reaching such a state, called stabilization, was studied  in \cite{S}
for 
a more general model on $\mathbb{R}^2$ with a finite number of initial seeds and arbitrary directions.

\begin{definition}[Escaping rays]
     Any ray $\mathcal{L}^{\pm}_{u,d_u}(N)$ in $G_{\mathcal{V}_N}(N)$  that reaches  the border of $\Lambda_N$  (i.e. $\mathcal{L}^{\pm}_{u,d_u}(N) \cap \partial \Lambda_N \neq \emptyset$)
     is called an {\normalfont escaping ray}.
\end{definition}

Thus, the total number of escaping rays quantifies the flow of the initial \enquote{activation} (as it was in the original model of neuronal activity) through a boundary. Recall that the number of initially active points is a random 
$\text{Poi}(\lambda N^2)$ variable. 
However, the number of escaping rays is much smaller than the average number of initially active points, as we shall establish next.

\begin{theorem}\label{th:escBound}
    Let $\mathcal{R}_N$ be the number of escaping
    rays in the final random quadrangulation $G_{\mathcal{V}_N}(N)$. There are constants $C_1, C_2$ such that
    \begin{equation*}
        C_1\sqrt{\lambda}N \leq \mathbb{E}(\mathcal{R}_N)< C_2\sqrt{\lambda}
        N
    \end{equation*}
    for any $\lambda > 0$ and for all large $N$.
\end{theorem}

\section{Proofs}\label{proofs}

\subsection{ Proof of Theorem \ref{T1}}
 
First we recall again the property of the Poisson process, which is of particular use here. 

\begin{remark}\label{Rem1}
    By the property of the Poisson process – independence for the processes over non-intersecting areas – for any finite set ${\cal{U}}$ of points in $\mathbb{R}^2$ the distribution of 
$\mathcal{V}$ 
conditional on $\cal{U} \in \mathcal{V}_N$ is equal to the distribution of the set $\mathcal{V} \cup \cal{U}$ (see, for example,  \cite{Kingman} or \cite{Cox}). 
Here we consider sets 
$\cal{U}$ consisting of one or two points only.
\end{remark}

From now on throughout the proof 
we can let $v\in \mathbb{R}^2$ and $d_v$ be fixed arbitrarily, and consider the dynamics of 
$$|\mathcal{L}^{+}_{v,d_v}(T)
\cup
\mathcal{L}^{-}_{v,d_v}(T)|=|\mathcal{L}^{+}_{v,d_v}(T)
|+|
\mathcal{L}^{-}_{v,d_v}(T)|={L}^{+}_{v,d_v}(T)
+{L}^{-}_{v,d_v}(T)
$$ in $G_{v \cup \mathcal{V}}(T)$, $T>0$, which is equivalent to carrying conditioning on $\{v\in \mathcal{V}\}$ throughout the proof. As we already derived in (\ref{lt})
for all $T>t>0$
\begin{equation}\label{R11}
 \mathbb{P}\{L^+_{v,d_v}(T)>t \} =
    \mathbb{P}\{L^+_{v,d_v}(t)=t \}.
\end{equation}

Taking into account the independence of $\mathcal{L}_{v,d_v}^+$ and $\mathcal{L}_{v,d_v}^-$ provided by Proposition \ref{Prop1}, the proof is carried out by proving the exponential decay for a ray on one side, $\mathcal{L}_{v,d_v}^+$, and is completed using an argument of symmetry.

Given 
a realization of $\mathcal{V}$, which assignment of directions $d_u$ for the  points  
\mbox{$u\in \mathbb{R}^2\setminus\{v\}$}
could yield 
the event $\{L_{v, d_v}^+(t)=t\}$? Note that $\{L_{v, d_v}^+(t)=t\}$ requires that
all rays from the points in $ C_{v, d_v}^+(t)$ (see Figure \ref{fig:coneDiamond}) which grow orthogonal towards $\mathcal{L}^+_{v, d_v}(t)$ are blocked by the lines 
 parallel to $\mathcal{L}^+_{v, d_v}(t)$.
For any Borel set $U\subseteq \mathbb{R}^2$ let us define the following event:
    \begin{equation}\label{RB}     
    {B(U):= \bigcap_{u\in U  \cap \mathcal{V}_N}\{ \mathcal{L}_{u, d_u}(t) \cap \mathcal{L}_{v, d_v}^+(t)=\emptyset\}}
    \end{equation}
    i.e., no rays from points in  $U$ reach 
$\mathcal{L}_{v, d_v}^+(t)$. Then we have a representation 
\begin{equation}\label{R16}
     \{{L}_{v, d_v}^+(t)=t\}=  
     B\left(C_{v, d_v}^+(t)\right)
    .
    \end{equation}

\begin{figure}[h]
        \centering
        \includegraphics[width=0.65\linewidth]{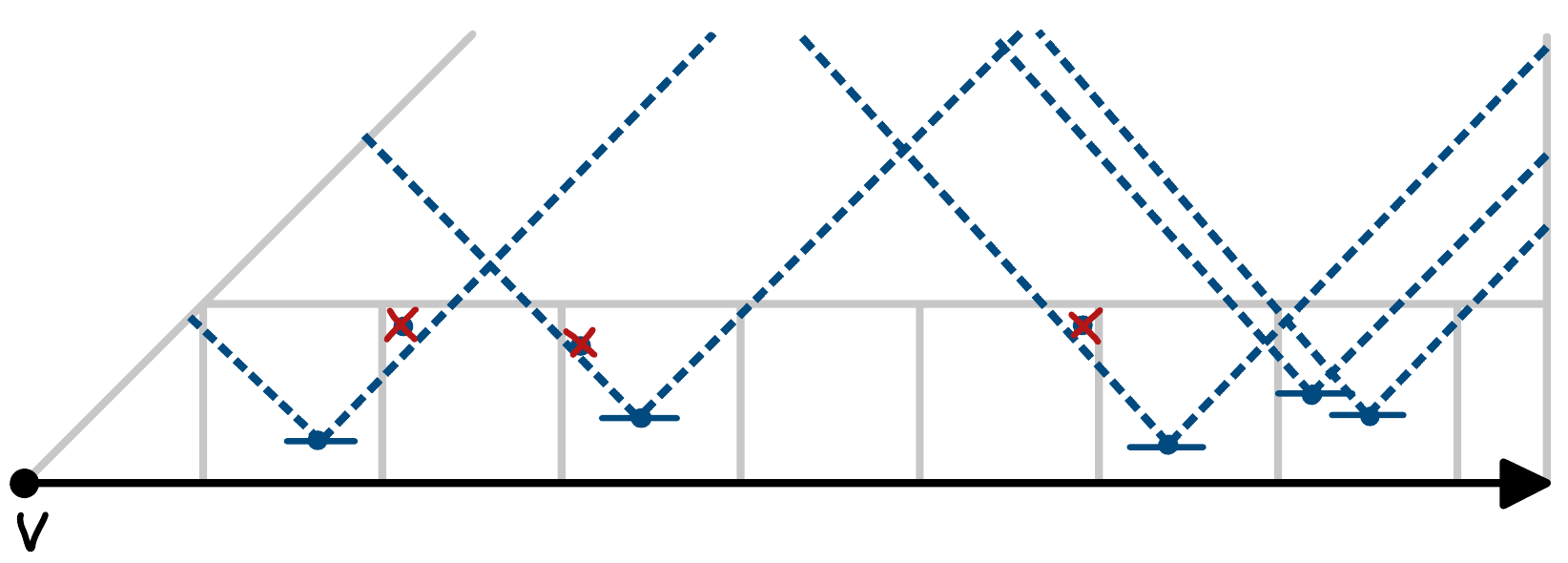}
        \caption{Illustration of \textit{tops}. The horizontal arrow of length $t$ represents  $\mathcal{L}_{v, d_v}^+(t)$. Blue points with small lines are the tops. Points that are crossed out in red are not tops, because the growth of lines from these points is blocked by the lines of other points. 
        } 
        \label{fig:IlluTopsBoxes}
    \end{figure}

To gain some intuition, consider Figure \ref{fig:IlluTopsBoxes} where we show only the left side along the ray; the right side would be similar. Figure \ref{fig:IlluTopsBoxes} illustrates that the lines parallel to the ray $\mathcal{L}_{v, d_v}^+(t)$  prevent orthogonal rays emanating from the points crossed out in red from reaching the ray $\mathcal{L}_{v, d_v}^+(t)$,
hence allowing event
$\{L_{v, d_v}^+(t)=t\}$ to happen.
  The dashed lines mark the borders of the areas within which the rays growing orthogonal to $d_v$ will not reach the ray $\mathcal{L}_{v, d_v}^+(t)$. This picture instructs us to identify the following pivotal positions in the configuration of points in a given set $\mathcal{V}$.

\begin{definition}[Top]\label{def:tops}
Let $v\in \mathbb{R}^2$ and a direction vector $d_v$ be fixed arbitrarily. For any given realization of $\mathcal{V}$ containing $v$, a point $u\in \mathcal{V}$
is called a {\normalfont top} with respect to $\mathcal{L}_{v, d_v}^\pm(t)$ in $\mathcal{V}
$,  if
\begin{equation}\label{R13}
\{L_{v, d_v}^\pm(t) = t\}\cap \{d_u \neq d_v\}= \emptyset 
\end{equation}
for all $t$ such that $\{L_{v, d_v}^\pm(t) = t\} \neq \emptyset$.
\end{definition} 

In other words, a point $u\in \mathcal{V}
$
is called a {\normalfont top} 
with respect to $\mathcal{L}_{v, d_v}^\pm(t)$
in $\mathcal{V}
$,
if, in order to satisfy condition $L_{v, d_v}^\pm(t)= t$,
the direction vector $d_u$ must  coincide with $d_v$.

Note that the condition (\ref{R13}) implies  that a top
with respect to $\mathcal{L}_{v, d_v}^\pm(t)$
could be found only within the area $C_{v, d_v}^\pm(t)$, because only the rays from the latter area might stop the growth of $\mathcal{L}_{v, d_v}^\pm$ before time $t$. 

To find bounds for the probability in (\ref{R11}) we shall first consider (given the set $\mathcal{V}$) the impact of points of $\mathcal{V}$ close to $\mathcal{L}_{v, d_v}^+(t)$.

\begin{figure}[htbp]
        \centering
        \includegraphics[width = 0.6 \linewidth]{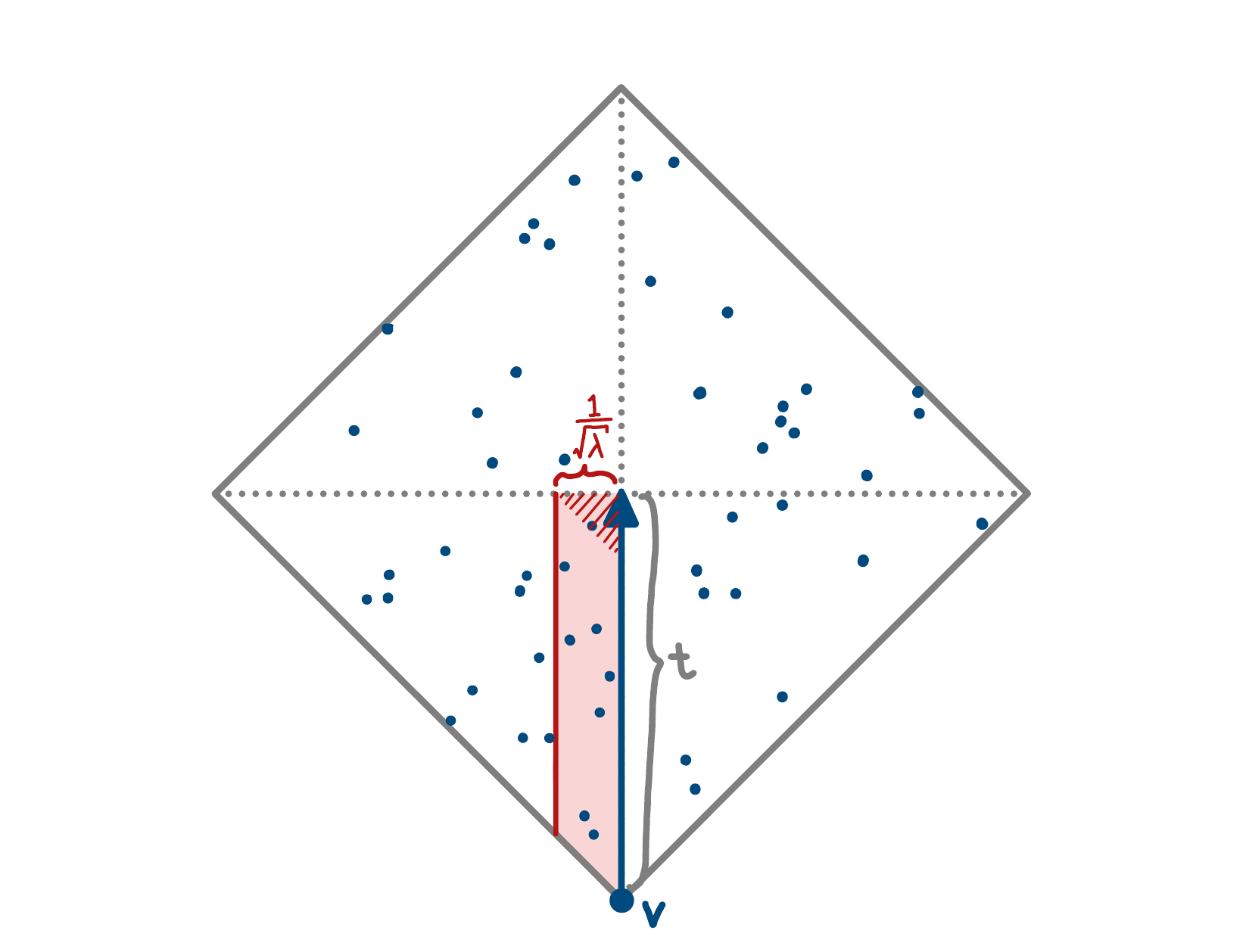}
        \caption{Illustration of the stripe $S_v^-(t)$ of width $\frac{1}{\sqrt{\lambda}}$, shaded in red. The (blue) dots represent points of $\mathcal{V}$.} 
        \label{fig:proofIllu}
    \end{figure}

    Define within $C_{v, d_v}^+(t)$ 
     areas $S_v^\pm(t)$ of width $\frac{1}{\sqrt{\lambda}}$ that are adjacent to $\mathcal{L}_{v, d_v}^+(t)$ on either side (one such area is shaded red in Figure \ref{fig:proofIllu}). 
     More precisely, $S_v^\pm(t)$
     are the following sets
     \begin{equation}\label{st}
         S_v^\pm(t):=\left\{u\in C_{v, d_v}^+(t): u\pm s\Bar{d}_v\in \mathcal{L}_{v, d_v}^+(t) \mbox{ for some } 0\leq s\leq \frac{1}{\sqrt{\lambda}}\right\},
     \end{equation}
     where $\Bar{d}_v$ denotes a vector orthogonal to ${d}_v$ and obtained by rotating ${d}_v$ by $+\frac{\pi}{2}$.
     We shall call sets $S_v^\pm(t)$ \textit{stripes}.
     
    By representation (\ref{R16}) we have
    \begin{equation}\label{R17}
     \{{L}_{v, d_v}^+(t)=t\}=
     B\left(C_{v, d_v}^+(t)\right)
    \subseteq 
    B(S_v^+(t) \cup S_v^-(t))
    =B(S_v^+(t)) \cap B(S_v^-(t)),
    \end{equation}
    which implies
\begin{equation}\label{R20}
        \mathbb{P}({L}_{v, d_v}^+(t)=t)\leq \mathbb{P}\left\{B(S_v^+(t))\cap B(S_v^-(t))\right\}.
    \end{equation} 
    
    Note that the events $B(S_v^+(t))$ and $B(S_v^-(t))$  are not independent due to the following observation.
    The rays, which are orthogonal to $d_v$ and originated at
    the points of 
    $\mathcal{V}_N$ in $S_v^+(t)$,
    do not reach $S_v^-(t)$ and vice versa. 
    However, the rays, which are parallel to $d_v$
    and originated in $S_v^\pm(t)$
    in a small square with the centre at point $v+td_v$ (i.e., the centre of $D_{v,d_v}^+(t)$) and the diagonal of the length
    $2/\sqrt{\lambda}$ along $d_v$,
    may indirectly affect the rays in 
    $S_v^\mp(t)$, which are orthogonal to $d_v$
  (see  the small striped area in the centre in Figure \ref{fig:proofIllu}). 

To avoid the above mentioned interference, let us cut off certain parts of $S_v^\pm(t)$ as follows. Within each stripe ${S}_v^\pm(t)$, define a rectangle
\begin{equation}\label{stt}
    \Tilde{S}_v^\pm(t) \subset {S}_v^\pm(t)
\end{equation}
of the size
\begin{equation}\label{size}
\frac{1}{\sqrt{\lambda}} \times \frac{\left[ {t}{\sqrt{\lambda}}\right]-2}{\sqrt{\lambda}},
\end{equation}
one side of which is the interval $[v+\frac{1}{\sqrt{\lambda}}d_v, v+\left(t-
\frac{1}{\sqrt{\lambda}}\right)d_v
]$.
Then the events 
$B\left( 
\Tilde{S}_v^+(t) 
\right)$ and $B
\left(\Tilde{S}_v^-(t) \right)$ are (conditionally) independent:
\begin{equation}\label{R18}
 \mathbb{P}\left\{B(\Tilde{S}_v^+(t))\cap B(\Tilde{S}_v^-(t))\right\}=\mathbb{P}\left\{B(\Tilde{S}_v^+(t))\right\}\mathbb{P}\left\{ B(\Tilde{S}_v^-(t))\right\},
 \end{equation}
and by (\ref{stt}) we have 
\begin{equation}\label{R17*}
    B\left(
    {S}_v^\pm(t)\right) \subseteq
    B\left(\Tilde{S}_v^\pm(t) \right).
\end{equation}
Now properties (\ref{R18}) and (\ref{R17*}) allow us to derive from (\ref{R20})
\begin{equation}\label{eq:finalArgUppExp}
        \mathbb{P}({L}_{v, d_v}^+(t)=t)\leq 
    \mathbb{P}\left\{B(\Tilde{S}_v^+(t))\right\}\mathbb{P}\left\{ B(\Tilde{S}_v^-(t))\right\}= \mathbb{P}^2\left\{ B(\Tilde{S}_v^-(t))\right\},
    \end{equation} 
    where the last equality is due to symmetry.

    Consider the last probability with the help of the notion \enquote{top} introduced above. 
    First, taking into account (\ref{size}) we shall divide the stripe $\Tilde{S}_v:= \Tilde{S}_v ^-(t)$ into 
    \begin{equation}\label{nt}
        n_t= \left[ {t}{\sqrt{\lambda}}\right]-2 \geq  {t}{\sqrt{\lambda}}-3
    \end{equation}
    boxes (squares)
    with side length $\frac{1}{\sqrt{\lambda}}$ as illustrated in Figure \ref{fig:IlluTopsBoxes}. Notice that the numbers of points of $\mathcal{V}$ in different boxes are independent and are distributed as $Poi(1)$.
    
    Let $W$ denote the (random) number of boxes in $\Tilde{S}_v$ that contain at least one point of $\mathcal{V}$, i.e. that are non-empty. Then there are at least $\max\left \{ 0 , \left( \left\lfloor \frac{W}{3}\right\rfloor-2\right)\right\} $ tops in area $\Tilde{S}_v$ because of the following arguments. 
    Note that the growth of a ray 
    in the direction orthogonal to $\mathcal{L}_{v, d_v}^+(t)$
    from any points in any box in $\Tilde{S}_v$ 
    can be stopped before it reaches 
    $\mathcal{L}_{v, d_v}^+(t)$
    by a line from another point only if the latter point belongs to either the same box or to one of the two neighbouring boxes. 
    Hence, 
     \begin{equation}\label{S5}
        M:=\#\{\text{tops in }\Tilde{S}_v\}
        \geq \left\lfloor \frac{W}{3}\right\rfloor-2.
    \end{equation}

Consider now the last probability in 
(\ref{eq:finalArgUppExp}). The event $B(\Tilde{S}_v^-(t))$ (recall (\ref{RB})) requires that all tops are assigned direction $d_v$, which has probability $\frac{1}{2}$ for each top independently.
    This gives us a bound
\begin{equation}\label{S4}
\begin{aligned}
    \mathbb{P}(B(\Tilde{S}_v^-(t))) &=\mathbb{E}\mathbb{P}(B(\Tilde{S}_v^-(t)) \mid M) \\
    &\leq\sum_{m=0}^K \mathbb{P}(M=m) + \sum_{m=K + 1}^{\infty}\mathbb{P}(B(\Tilde{S}_v^-(t))|M=m)\\
&\leq \mathbb{P}(M\leq K) + \sum_{m=K + 1}^{\infty} \left(\frac{1}{2}\right)^m \\
&= \mathbb{P}(M\leq K) +\left(\frac{1}{2}\right)^K
\end{aligned}
\end{equation}
for an arbitrary $K>0$, which will be chosen in some optimal way. Observe that some probabilities here and below are unconditional; see Remark \ref{Rem1}.
    By (\ref{S5})
\begin{equation}\notag
        \mathbb{P}(M\leq K) \leq 
        \mathbb{P}(W\leq 3\gamma +8),
    \end{equation}
    where $W$ is the number of non-empty boxes among $n_t$ boxes. Recall that the area of each box is $1/\lambda$, and the points are distributed as a Poisson process with intensity $\lambda$, implying that the numbers of points in different boxes are independent and distributed as $\text{Poi}(1)$. Hence, the probability that a box is non-empty is $1-e^{-1}$, and the total number of non-empty boxes has a Binomial distribution:
\begin{equation}\label{S7*}
    W \sim \text{Bin}(n_t, 1-e^{-1}).
\end{equation}
By the Chernoff bound for all $0<\vartheta<1$
\begin{equation}\notag
   \mathbb{P} \{ W \leq \mathbb{E}W (1- \vartheta)\}\leq \exp\left\{-\frac{\vartheta ^2}{2}\mathbb{E}W\right\}.
\end{equation}
Choosing now $K$ in (\ref{S4}) so that 
\begin{equation}\notag
3K +8 = \left[ \mathbb{E}W (1- \vartheta)\right]
\end{equation}
for $0<\vartheta<1$ to be chosen later, we derive from (\ref{S4})
\begin{equation}\label{S7}
     \mathbb{P}(B(\Tilde{S}_v^-(t))) \leq
     \exp\left(-\frac{1}{2}\vartheta^2
     \mathbb{E}W\right) 
     + \exp\left(-\frac{\ln(2)}{3}(1-\vartheta)
     \mathbb{E}W+3\ln(2)\right),
\end{equation}
where
\begin{equation}\notag
    \mathbb{E}W = n_t(1-e^{-1})=
    (\lfloor t\sqrt{\lambda}\rfloor - 2)(1-e^{-1})\geq (t\sqrt{\lambda} - 3)(1-e^{-1}).
    \end{equation}
    As we are interested in the large values for $t$, the optimal value for $\vartheta$
matching the exponents in  (\ref{S7}) must satisfy
    \begin{equation}\notag
        \frac{1}{2}\vartheta^2=\frac{\ln(2)}{3}(1-\vartheta),
    \end{equation}
    which is 
    \begin{equation}\notag
        \vartheta =-\frac{\ln(2)}{3}+\sqrt{\frac{\ln(2)^2}{9}+\frac{2\ln(2)}{   3}}.
    \end{equation} 
    With the last choice of $\vartheta$ the bound in (\ref{S7}) reads as 
    \begin{equation}\notag
        \mathbb{P}(B(\Tilde{S}_v^-(t)))
        \end{equation}
        \[
        \leq\left(\exp\left(\frac{3}{2}(1-e^{-1})\vartheta^2\right) +  \left(\frac{1}{2}\right)^{-(3+(1-\vartheta )(1-\frac{1}{e}))} \right) \exp
        \left(-
        \frac{\vartheta^2}{2}(1-e^{-1})\sqrt{\lambda}t\right)
    =:ce^{-\alpha \sqrt{\lambda}t}\]
    where $c$ and $\alpha$ are some positive constants independent of $\lambda$.
Substituting the last bound into (\ref{eq:finalArgUppExp}) 
we get
\begin{equation} \notag
        \mathbb{P}({L}_{v, d_v}^+(t)=t)\leq c^2e^{-2\alpha \sqrt{\lambda}t}.
    \end{equation}
    
    By symmetry, the same bound holds for  $L_{v, d_v}^-(t)$ as well:
    \[ \mathbb{P}({L}_{v, d_v}^-(t)=t)\leq c^2e^{-2\alpha \sqrt{\lambda}t}.\]
    Hence, for all $t<T$
    \begin{equation} \label{S20}
   \mathbb{P}\{
      |\mathcal{L}^{+}_{v,d_v}(T)
\cup
\mathcal{L}^{-}_{v,d_v}(T)|
      >t \mid v\in\mathcal{V}\}
     \end{equation}
     \[\leq 
     \mathbb{P}\{
      |\mathcal{L}^{+}_{v,d_v}(T)|>t/2
\mid v\in\mathcal{V}\}+
\mathbb{P}\{
      |\mathcal{L}^{-}_{v,d_v}(T)|>t/2
\mid v\in\mathcal{V}\}
    \leq 2c^2e^{-\alpha \sqrt{\lambda}t}
    ,\]
    which is the upper bound in the statement of the theorem.

    Let us now turn to the lower bound.

Recall the formal representation (\ref{R16}). This time we shall consider potential rays from the entire area of the cone $C_{v, d_v}^+(t)\setminus [v, v+td_v]$ and not only from the stripes as in the previous case.
Let us
split the cone into the left and right parts with respect to the direction of $d_v$ (the middle ray may go to either part) and denote the resulting triangles $T_{v}^{\pm}(t)$, respectively:
\[C_{v, d_v}^+(t) =:  T_{v}^-(t) \cup T_{v}^+(t).\]
For the stripes defined in (\ref{st}), it holds that $S_{v}^\pm(t)\subset T_{v}^\pm(t)
$.
Then, representation (\ref{R16})
yields
 \begin{equation}\label{2B}
        \mathbb{P}({L}_{v, d_v}^+(t)=t)= \mathbb{P}\left\{B(T_v^+(t))\cap B(T_v^-(t))\right\}.
    \end{equation}
To get a lower bound for this probability, we shall find events which yield, correspondingly, each of the events on the right side of (\ref{2B}).

\begin{figure}[h]
    \centering
    \includegraphics[width =0.5\linewidth]{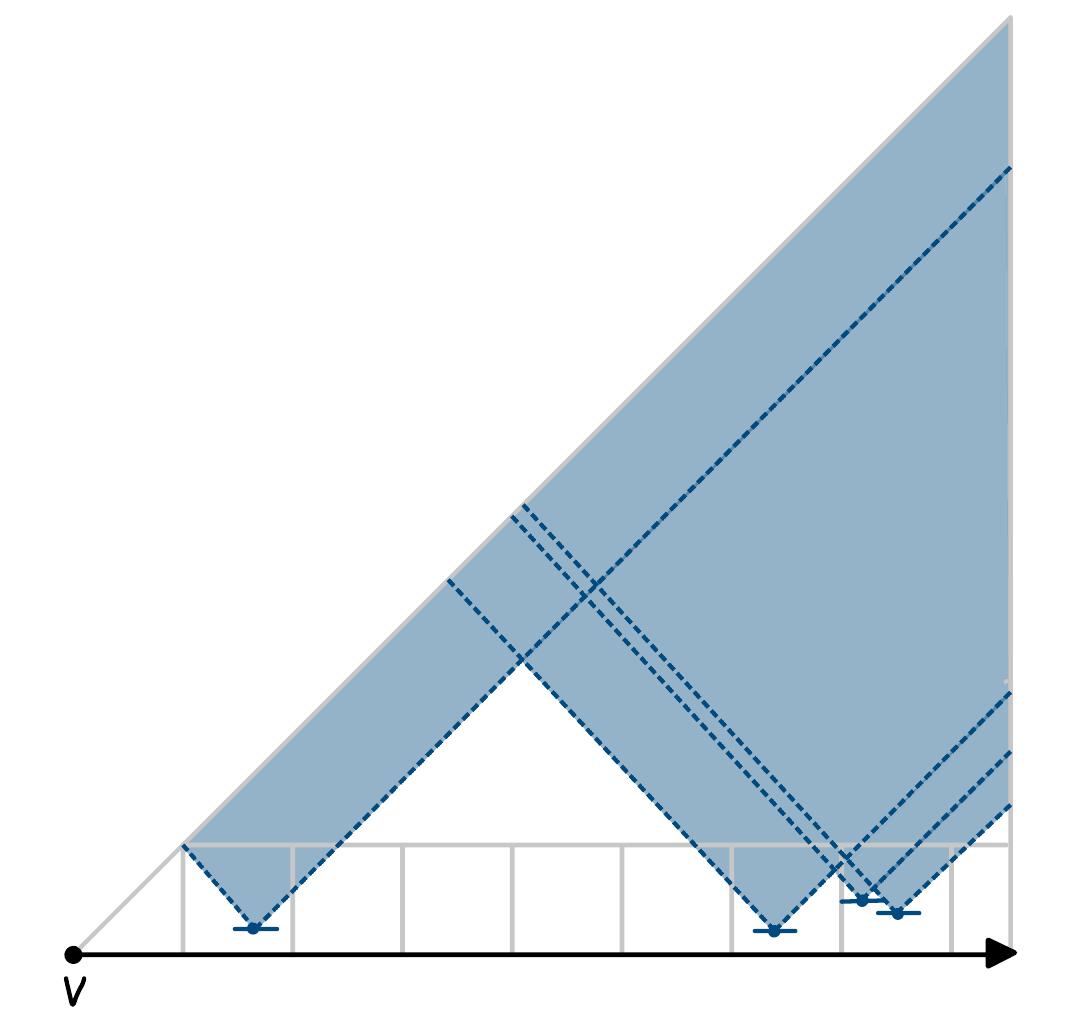}
    \caption{The ray is $\mathcal{L}_{v,d_v}^+(t)$. The entire triangle above the ray is $T_v^-(t)$. The points with lines are the tops. The area where the directions assigned to the points do not affect $\{L_{v,d_v}(t)=t\}$ is shaded in blue.} 
    \label{fig:IlluLowerDecayBasic}
\end{figure}

We use again a rectangle within a stripe on one side of the ray $\mathcal{L}_{v,d_v}^+(t)$ 
defined in (\ref{stt}), say:
\begin{equation}\label{stt*}
    \Tilde{S}_v^-(t) \subset {S}_v^-(t) \subset {T}_v^-(t);
\end{equation} 
 see Figure \ref{fig:IlluLowerDecayBasic}. 
 
If a given configuration $\mathcal{V}$ 
has some points in $\Tilde{S}_v^-(t)$, then some of those are the \enquote{tops} with respect to $\mathcal{L}_{v,d_v}^+(t)$ (see Definition \ref{def:tops}). 
For a shorthand notation, let us denote here for a given $\mathcal{V}$ and a set $A\subseteq D_{v,d_v}^+(t)$ 
\begin{equation}\label{Topnot}
    {\cal{T}}(A):=\{\mbox{the tops in $A$ 
 with respect to } \mathcal{L}_{v,d_v}^+(t)\}
\end{equation}
 the set of the tops in $A$ 
 with respect to the ray $\mathcal{L}_{v,d_v}^+(t)$.

Observe that if all the tops in ${\cal {T}}\left(\Tilde{S}_{v}^-(t)\right)$ are assigned the same vector $d_v$, i.e., 
\begin{equation}\label{dTop}
    d_u=d_v, \ \ u\in {\cal {T}}\left(\Tilde{S}_{v}^-(t)\right),
\end{equation}
then the lines from these tops will prevent crossing of $\mathcal{L}_{v,d_v}^+(t)$ 
by any ray
 from their dependence area,
 which in notation of Definition \ref{depA} is
\begin{equation}\label{U}
      \bigcup_{u\in  {\cal {T}}\left(\Tilde{S}_{v}^-(t)\right)}D^+_{u,\Bar{d}_v}(t)
 \end{equation}
 (see the shaded area above the tops in Figure \ref{fig:IlluLowerDecayBasic}). Hence, assuming further that all  other lines outside of the area (\ref{U}), i.e., in 
 \begin{equation}\label{UX}
     U_v^-:=T_{v}^-(t) \setminus
     \left( \bigcup_{u\in  {\cal {T}}\left(\Tilde{S}_{v}^-(t)\right)}D^+_{u,\Bar{d}_v}(t)\right)
 \end{equation}
 (the white area above the ray in the Figure \ref{fig:IlluLowerDecayBasic}),
 are parallel to $d_v$, and thus  neither cross $\mathcal{L}_{v,d_v}^+(t)$,
this will 
yield the event $B(T_v^-(t))$.

We collect all these assumptions  and define the following event which guarantees that no line from $T_v^-(t)$ crosses $\mathcal{L}_{v,d_v}^+(t)$:
\begin{equation}\label{B-}
    B^-:=\left\{
    d_u=d_v, \mbox{ for all }
    u\in U_v^- \cup {\cal{T}}
    \left(\Tilde{S}_{v}^-(t)\right)
    \right\}\subseteq B(T_v^-(t)).
\end{equation}
Repeating the same construction, we define a similar event for the other area $T_v^+(t)$ of the ray $\mathcal{L}_{v,d_v}^+(t)$:
\begin{equation}\label{B+}
    B^+\subseteq B(T_v^+(t)).
\end{equation}
By the symmetry and construction of $\Tilde{S}_{v}^\pm(t)$, the events $B^\pm$ are independent and have the same distribution.
This gives us the following lower bound for the probability \mbox{in (\ref{2B})}
\begin{equation}\label{R30}
        \mathbb{P}\left\{B(T_v^+(t))\cap B(T_v^-(t))\right\}
        \geq \mathbb{P}\left\{B^+\cap B^-\right\}= 
      \mathbb{P}\left\{B^+\right\}  \mathbb{P}\left\{ B^-\right\}
        .
    \end{equation}
    
We shall derive now the lower bound for 
$\mathbb{P}\left\{B^-\right\}$,
where $B^-$
 is defined by
 (\ref{B-}). The idea is the following. Let $K_t$ denote the number of tops in ${\cal{T}}
    \left(\Tilde{S}_{v}^-(t)\right)$ and let 
    \begin{equation}\label{Y}
    Y:=\left|U_v^- \cap \mathcal{V}\right|
    \end{equation}
    denote the number of points of $\mathcal{V}$ in $U_v^-$ (defined in (\ref{UX})).
Then by definition (\ref{B-})
\begin{equation}\label{R31}
    \mathbb{P}\left\{B^-
    \right\}=\mathbb{E}
    \mathbb{P}
    \left\{B^-\mid Y, K_t
    \right\}= \mathbb{E}\left(\frac{1}{2}\right)^{Y+K_t},
\end{equation}
which simply takes into account that the points are assigned direction $d_v$ independently and with probability $1/2$. Using Jensen's inequality we get from (\ref{R31}) 
\begin{equation}\label{R41}
    \mathbb{P}\left\{B^-
    \right\}\geq  2^{-\mathbb{E}(Y+K_t)}.
\end{equation}
By definition (\ref{Y}) and by the definition of the Poisson process, the distribution of $Y$ given the area $A=A(U_v^-)$ of the set $U_v^-$ follows the Poisson distribution:
\begin{equation}\label{Ydist}
    Y \sim \text{Poi}(\lambda A),
\end{equation}
hence,
\begin{equation}\label{R42}
    \mathbb{E}Y=\lambda \mathbb{E}A(U_v^-).
\end{equation}
To bound the expectation 
in (\ref{R42})
we shall use the following construction.
We divide again the set 
 $\Tilde{S}_v^-(t)$ into $n_t$ boxes (see (\ref{nt})), and 
set $I=(I_1, \ldots , I_{n_t})$,  where a random variable $I_i$ denotes the indicator that the $i$-th
box is empty. By the properties of the Poisson process, $I_i$ are $i.i.d.$ Bernoulli random variables with a parameter $1/e$.

\begin{figure}[h]
    \centering
    \includegraphics[width =0.75\linewidth]{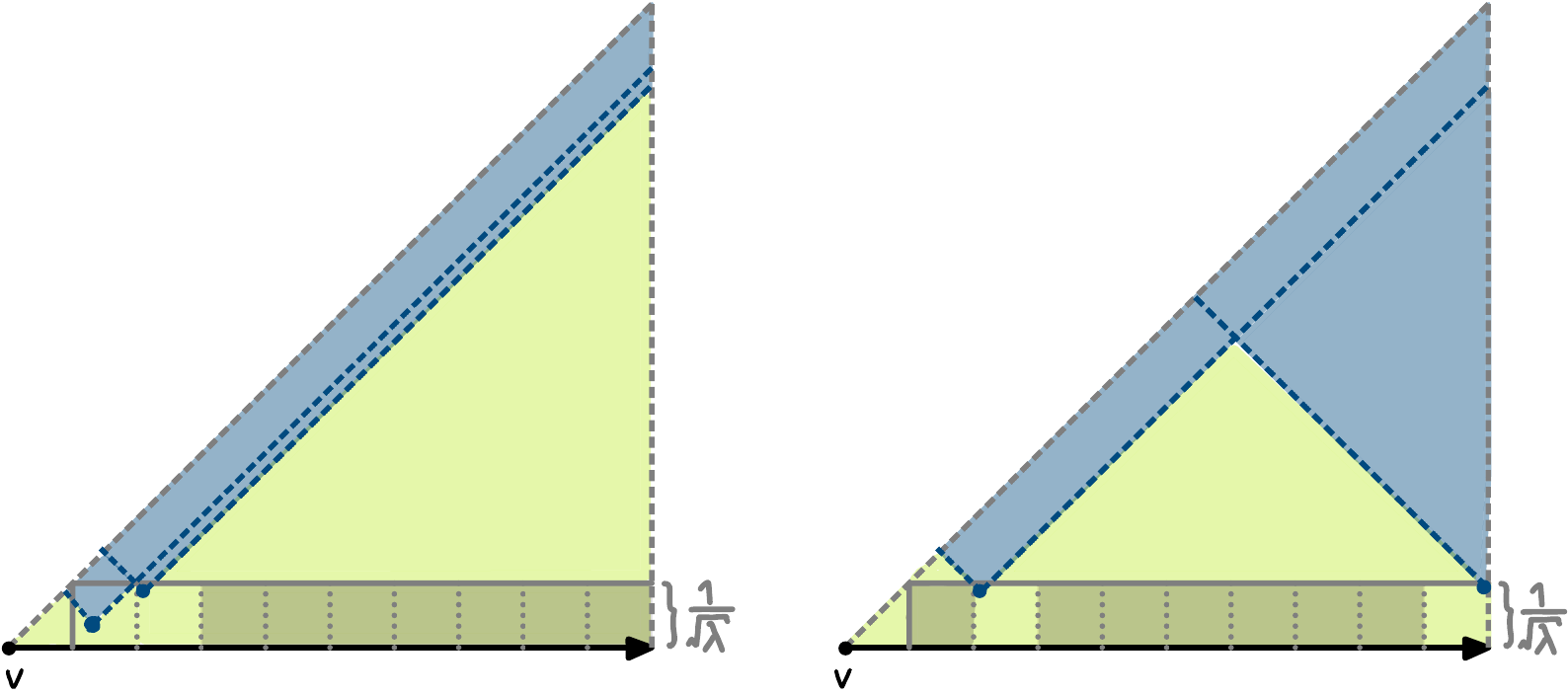}
    \caption{$U_v^-$ shaded in light green, with the example of $7$ empty boxes (shaded darker) in different positions.} 
    \label{fig:IlluLowerDecay}
\end{figure}

Given  that $I$ contains at least one zero, i.e.,  $I\neq (1, \ldots, 1)$,  let $X_1$, $X_2$, $\ldots, X_T$ be  the lengths of all series of consecutive zeros in $I$, where $T=T(I)$ is the number of such series of zeros (or empty boxes) separated by values $1$ (non-empty boxes); see some examples in Figure \ref{fig:IlluLowerDecay}. Observe that by the Strong Law of Large numbers, as $n_t\rightarrow \infty$
\begin{equation}\label{asconv}
    \frac{n_t}{T(I)} \stackrel{a.s.}{\rightarrow} \mathbb{E}(\xi _0+\xi _1),
\end{equation}
where $\xi_0 $ has the First success distribution Fs$\left(\frac{1}{e}\right)$, i.e.,
\[\mathbb{P}\{\xi _0=k\}=\left(1-\frac{1}{e}\right)^{k-1}\frac{1}{e}, \ \ k=1,2, \ldots , \]
and $\xi_1 $ correspondingly has the First success distribution Fs$\left(1-\frac{1}{e}\right)$. The last two distributions are the distributions of the number of consecutive zeros and the number of consecutive ones in the unbounded sequence of independent $I_k, k\geq 1$ (when $n_t \rightarrow \infty$).
    
Figure \ref{fig:IlluLowerDecay} helps us to realize, that the 
area $U_v^-(t)$ is represented by the triangles over the intervals of consecutive zeros (stretched at most by $1$ unit into both sides); thus, given the vector $I$, we
derive the following bound 
\begin{equation}\label{eq:areaSvw}
        A(U_v^-(t)) \leq \sum_{i=1}^{T(I)}
        \frac{1}{2\lambda}(X_i+2)^2 +\frac{2n_t}{\lambda},
    \end{equation} 
    where the last term 
    corresponds to the case when there are no empty boxes: it 
    is the upper bound for the area (of a strip of width $1/\sqrt{\lambda}$ along
    $\widetilde{S}_v^-(t)$)
    not covered by $U^-_v(t)$.
    Bound (\ref{eq:areaSvw}) together with observation (\ref{asconv}) implies
\begin{equation}\label{R37*}
  \mathbb{E}Y= \lambda \mathbb{E}A(U_v^-(t)) \leq C_1 {n_t},
\end{equation}
where $C_1$ is some positive constant independent of $\lambda$, and together with 
(\ref{nt}) this yields
\begin{equation}\label{R37}
  \mathbb{E}Y \leq C_1 {\sqrt{\lambda}}t.
\end{equation}

To bound $\mathbb{E}K_t$ we observe that $K_t$ is at most the number of all points  of $\mathcal{V}$ in the set
$\Tilde{S}_v^-(t)$, which is again the Poisson random variable
Poi$(\lambda A(\Tilde{S}_v^-(t)))$, where
$A(\Tilde{S}_v^-(t))$
denotes the area of $\Tilde{S}_v^-(t)$. Thus we have (recalling the construction around (\ref{size}))
\begin{equation}\label{R38}
    \mathbb{E}K_t = {\lambda}A(\Tilde{S}_v^-(t))\leq {\sqrt{\lambda}}t.
\end{equation}
Substituting the bounds (\ref{R37}) and  (\ref{R38}) into 
(\ref{R41}) and recalling definition (\ref{nt}), we get 
\begin{equation}\label{R42*}
    \mathbb{P}\left\{B^-
    \right\}\geq   c_1 e^{-\alpha_1 \sqrt{\lambda}t}
\end{equation}
for some constants
$c_1 ,  \alpha _1$.

By symmetry, the same bound (\ref{R42*}) holds for $\mathbb{P}\left\{B^+
    \right\}$, and hence the last bound allows us to derive from (\ref{R30})
\begin{equation}\label{R43}
        \mathbb{P}\left\{B(T_v^+(t))\cap B(T_v^-(t))\right\}
        \geq c_1^2 e^{-2 \alpha_1 \sqrt{\lambda}t}.
    \end{equation}
    Making use of this result in (\ref{2B}) we get
 \begin{equation}\label{S15}
        \mathbb{P}({L}_{v, d_v}^+(t)=t)= \mathbb{P}\left\{B(T_v^+(t))\cap B(T_v^-(t))\right\}
    \geq c_1^2 e^{-2 \alpha_1 \sqrt{\lambda}t}    .
    \end{equation}
    
By symmetry, the same bound holds when ${L}_{v, d_v}^+(t)$ is replaced by ${L}_{v, d_v}^-(t)$.
Therefore
\begin{equation}\notag
    \mathbb{P} \{{L}_{v, d_v}^+(N)+{L}_{v, d_v}^-(N)>t\} \geq c_1^2e^{-2\alpha_1 \sqrt{\lambda}t}.
\end{equation}
This together with the upper bound (\ref{S20}) completes the proof of Theorem \ref{T1}. 
\hfill$\Box$

\subsection{Proof of Proposition \ref{Prop2}}

Following Remark \ref{Rem1} we may fix $u,v$ and consider the growth of rays from these points in the environment of a Poisson point process.

\begin{figure}[h]
    \centering
    \includegraphics[width=0.5\linewidth]{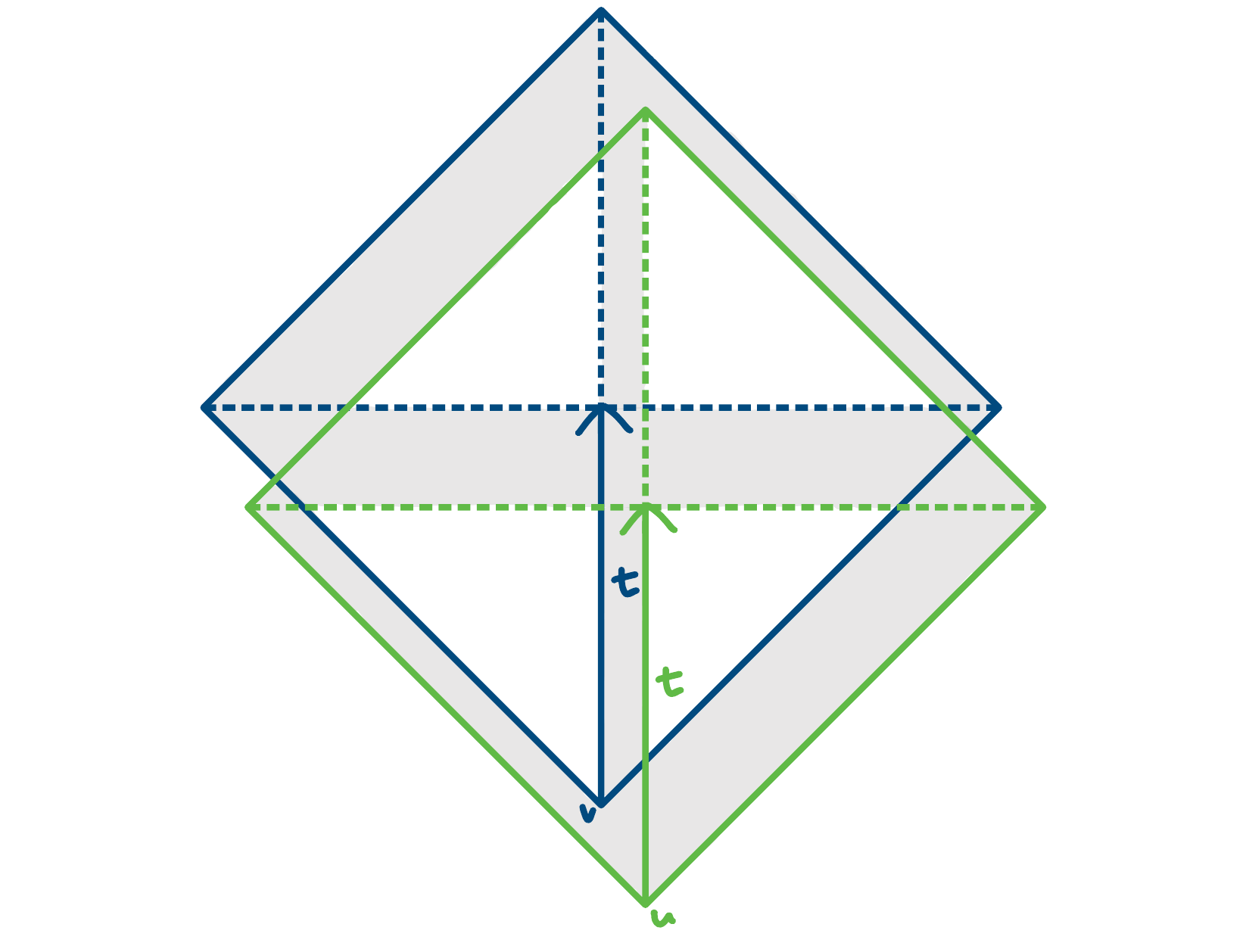}
    \caption{Illustration of the constructed area $\mathcal{B}$, shaded in grey.}
    \label{fig:IlluPosCov}
\end{figure}

For any $u,v$ with $\|u-v\|_1=l$ and 
with direction $d_u=d_v=d$ fixed arbitrarily,
consider the areas of dependence $$ D_{v,d}^+(t), \ D_{u,d}^+(t), \ t>0.$$ 
Let $\mathcal{C}$ denote
 the central area between the rays $\mathcal{L}_{v,d}^+$ and $\mathcal{L}_{u,d}^+$ within the union \mbox{$D_{v,d}^+(t) \cup D_{u,d}^+(t)$} (see Figure \ref{fig:IlluPosCov}), and then define another set
\[\mathcal{B}=\mathcal{C} \cup \left( D_{v,d}^+(t) \triangle D_{u,d}^+(t)\right)\cup \left( C_v^+(t) \triangle C_u^+(t)\right) . \]
This is illustrated in Figure \ref{fig:IlluPosCov}, where $\mathcal{B}$ is shaded in grey. Then we observe that (recall definition (\ref{depA}))
\begin{equation}\label{S24}
A_{v,d}^+(t)
\cap \{\mathcal{B} \cap \mathcal{V} = \emptyset\}=A_{u,d}^+(t)
\cap \{\mathcal{B} \cap \mathcal{V} = \emptyset\}
\end{equation}
\[=
A_{v,d}^+(t)
\cap A_{u,d}^+(t)
\cap \{\mathcal{B} \cap \mathcal{V} = \emptyset\}
\subseteq 
A_{v,d}^+(t)
\cap A_{u,d}^+(t),
\]
which gives us
\begin{equation}\label{S22N}
  \mathbb{P}  \left\{
  A_{v,d}^+(t)
\cap A_{u,d}^+(t)
  \right\}
\geq 
\mathbb{P} \left\{ 
A_{v,d}^+(t)
\cap \{\mathcal{B} \cap \mathcal{V} = \emptyset\}
\right\}.
\end{equation}

The number of points in $\{\mathcal{B} \cap \mathcal{V}\}$ has a Poisson distribution; hence
\begin{equation}\label{S25*}
  \mathbb{P} \left\{ \mathcal{B} \cap \mathcal{V} = \emptyset\right\}= e^{-\lambda Area(\mathcal{B} )},
\end{equation}
where $Area(\mathcal{B} )$  denotes the area of $\mathcal{B} $. The set $\mathcal{B} $ goes along the perimeter of $D_v^+(t)$ and its two diagonals; this yields a bound
\begin{equation*}
    Area(\mathcal{B} )<4(1+\sqrt{2})(t+l)l=:c(t+l)l,
\end{equation*} where $l$ is the distance between $u$ and $v$.
Thus, together with (\ref{S25*}), this 
gives us
\begin{equation}\label{S25}
  \mathbb{P} \left\{ \mathcal{B} \cap \mathcal{V} = \emptyset\right\}\geq e^{-\lambda c(t+l)l}.
\end{equation}
With the help of (\ref{S25}), we proceed with computing  the last probability in (\ref{S22N}):
\begin{equation}\label{64}
\mathbb{P}  \left\{
  A_{v,d}^+(t)
\cap A_{u,d}^+(t)
  \right\}
\geq \mathbb{P} \left\{ A_{v,d}^+(t)
\cap \{\mathcal{B} \cap \mathcal{V}_N = \emptyset\} \right\}
\end{equation}
\[=
\mathbb{P} \left\{ A_{v,d}^+(t) \right\}-
\mathbb{P} \left\{ A_{v,d}^+(t)
\cap \{\mathcal{B} \cap \mathcal{V}_N \neq \emptyset\}
\right\}\]
\[
\geq
\mathbb{P} \left\{ A_{v,d}^+(t)\right\}-
\mathbb{P} \left\{\mathcal{B} \cap  \mathcal{V}_N \neq \emptyset\right\}
\geq
\mathbb{P} \left\{ A_{v,d}^+(t) \right\}-\left(1-e^{-c\lambda (t+l)l}\right).
\]

We are ready now to bound the covariance (see definition (\ref{CABu})) as follows
\begin{equation}\label{Cov1R}
    {\bf C}\left(A^+_{u,d_u}(t), A^+_{v,d_v}(t)\right)=
    \mathbb{P}\{A^+_{u,d_u}(t)\cap A^+_{v,d_v}(t)\}-\mathbb{P}\{A^+_{u,d_u}(t)\}\mathbb{P}\{A^+_{v,d_v}(t)\}
\end{equation}
\[\geq
\mathbb{P} \left\{ A^+_{u,d_u}(t)\right\}(1-\mathbb{P}\{A^+_{v,d_v}(t)\})-\left(1-e^{-c\lambda (t+l)l}\right).
\]
Substituting the upper and the lower bounds established in Theorem \ref{T1} in the last formula, we obtain
\begin{equation}\notag
    {\bf C}\left(A^+_{u,d_u}(t), A^+_{v,d_v}(t)\right)
    \end{equation}
\[\geq c_1 e^{-\alpha_1 {\sqrt{\lambda}}t}
(1-c_2 e^{-\alpha_2 {\sqrt{\lambda}}t})-c\lambda (t+l)l,
\]
which is positive for all large $t$ and
\[l<c_3\frac{e^{-\alpha_1 {\sqrt{\lambda}}t}}{\lambda t},\]
where $c_3>0$ is a positive constant uniform in $\lambda$. \hfill$\Box$

\subsection{Proof of Theorem \ref{th:covExpBound}}
Let $u$ and $v$ together with the direction vectors $d_u$ and $d_v$ be fixed arbitrarily. Conditionally on $v,u \in \mathcal{V}$ (see also Remark \ref{Rem1}),
consider  events $A_{u,d_u}^{\pm}(t)$

Assume that $t>\frac{1}{4}\|u-v\|_1$. Otherwise, the statement is trivial due to \mbox{Corollary \ref{le:distForInd}}.
Note that for any $t'<t$ and $u \in \mathcal{V}$, we have
\begin{equation}\label{S21}
A_{u,d_u}^{\pm}(t) \subseteq A_{u,d_u}^{\pm}(t'),
\end{equation}
which implies
\begin{equation}\label{S22}
\mathbb{P}\left\{
A_{u,d_u}^{\pm}(t) \cap A_{v,d_v}^{\pm}(t) \mid v,u \in \mathcal{V}\right\} 
 \leq \mathbb{P}\left\{A_{u,d_u}^{\pm}(t') \cap A_{v,d_v}^{\pm}(t')\mid v,u \in \mathcal{V}\right\}.
 \end{equation}
Therefore, setting now
    \begin{equation*}
        {t}'=\frac{1}{4}\|u-v\|_1,
    \end{equation*}
and
    applying Corollary \ref{le:distForInd} on the independence of the events on the right-hand side in (\ref{S22}) we continue (\ref{S22}) as follows
\begin{equation}\label{S22*}
\mathbb{P}\left\{
A_{u,d_u}^{\pm}(t) \cap A_{v,d_v}^{\pm}(t) \mid v,u \in \mathcal{V}\right\} 
 \leq \mathbb{P}\left\{A_{u,d_u}^{\pm}(t') \cap A_{v,d_v}^{\pm}(t')\mid v,u \in \mathcal{V}\right\} 
 \end{equation}
 \[=\mathbb{P}\left\{A_{u,d_u}^{\pm}(t') \mid u \in \mathcal{V}\right\}
 \mathbb{P}\left\{
  A_{v,d_v}^{\pm}(t')\mid v \in \mathcal{V}\right\}.
\]
Hence, (\ref{S22*}) yields
\begin{equation}\label{eq:covDec1}
      \left|
      {\bf C}\left(A_{u,d_u}^{\pm}(t), A_{v,d_v}^{\pm}(t)\mid u,v  \in\mathcal{V} \right)
      \right|
    \end{equation}
\[ \leq 
\mathbb{P}\left\{A_{u,d_u}^{\pm}(t') \mid u \in \mathcal{V}\right\}
 \mathbb{P}\left\{
  A_{v,d_v}^{\pm}(t')\mid v \in \mathcal{V}\right\}\]
  \[+\mathbb{P}\left\{A_{u,d_u}^{\pm}(t) \mid u,v \in \mathcal{V}\right\}
 \mathbb{P}\left\{
  A_{v,d_v}^{\pm}(t)\mid u,v \in \mathcal{V}\right\},
\]
which, with the help of (\ref{S21}) again, gives us
\begin{equation}\label{eq:covDec2}
   \left|{\bf C}\left(A_{u,d_u}^{\pm}(t), A_{v,d_v}^{\pm}(t)\mid u,v  \in\mathcal{V} \right)\right|\leq  2 \left(\mathbb{P}\left\{A_{u,d_u}^{\pm}(t') \mid u \in \mathcal{V}\right\}\right)^2,
    \end{equation}
where we also made use of the symmetry properties in the model.
    Finally, applying Theorem \ref{T1} we derive from (\ref{eq:covDec2})
\begin{equation}\notag
     \left|{\bf C}\left(A_{u,d_u}^{\pm}(t), A_{v,d_v}^{\pm}(t)\mid u,v  \in\mathcal{V} \right)\right|\leq 2c_2^2 e^{-{\alpha}_2{\sqrt{\lambda}}\frac{1}{2}\|u-v\|_1},
    \end{equation}
    which is the statement of the theorem. \hfill$\Box$

\subsection{Proof of Theorem \ref{th:escBound}}

Given a set $\mathcal{V}_N$, for any  $v\in \mathcal{V}_N$, let us introduce the events
\[B_v^{\pm}:=\{ \mathcal{L}_{v,d_v}^{\pm}(N) \mbox{ reaches the border of } \Lambda_N\}. \]
Our aim is to  study the expectation of
\begin{equation}\notag
    \mathcal{R}_N:= \sum_{v\in \mathcal{V}_N}
    \left(\mathbbm{1}\{B_v^+\} + \mathbbm{1}\{B_v^-\} \right).
\end{equation}
Given a set $\mathcal{V}_N$, let  $U$ denote a random point chosen uniformly in $\mathcal{V}_N$. Then we can write
\begin{equation}\label{ERn}
   \mathbbm{E} \mathcal{R}_N = \mathbbm{E} \sum_{v\in \mathcal{V}_N}
    \left(\mathbbm{1}\{B_v^+\} + \mathbbm{1}\{B_v^-\} \right)
    \end{equation}
    \[
    = \mathbbm{E}\left( |\mathcal{V}_N|\mathbbm{E}\left\{\left.
    \left(\mathbbm{1}\{B_U^+\} + \mathbbm{1}\{B_U^-\} \right)\right||\mathcal{V}_N|\right\}\right).
\]
Consider events $B_U^{\pm}$. Letting $\eta_U^{\pm}$ denote the distance from $U$ to the boundary of $\Lambda_{N}$ in the direction of $\mathcal{L}_{U,d_U}^{\pm}$ correspondingly, we have
\begin{equation}\notag
    B_U^{\pm}=\{{L}_{U, d_U}^{\pm}(N)\geq \eta_U^{\pm}\}.
\end{equation}
Hence, we can rewrite 
(\ref{ERn}) as 
\begin{equation}\label{bd2}
\mathbbm{E} \mathcal{R}_N    
\end{equation}
\[=   
\mathbbm{E}\left( |\mathcal{V}_N|\mathbbm{P}\left\{\left.
    {L}_{U, d_U}^{+}(N)\geq \eta_U^{+}
    \right||\mathcal{V}_N|\right\}\right)
    +\mathbbm{E}\left( |\mathcal{V}_N|\mathbbm{P}\left\{\left.
    {L}_{U, d_U}^{-}(N)\geq \eta_U^{-}
    \right||\mathcal{V}_N|\right\}\right).\]

We shall make use of the property of the Poisson process that
conditionally on a positive number $|V_N|$ the distance $\eta_U^+$
is uniformly distributed on the entire interval $[0,N]$. 
This allows us to rewrite
\begin{equation}\label{bd5}  
\mathbbm{E}\left( |\mathcal{V}_N|\mathbbm{P}\left\{\left.
    {L}_{U, d_U}^{+}(N)\geq \eta_U^{+}
    \right||\mathcal{V}_N|\right\}\right)
\end{equation}
    \[=
\int_0^N \frac{1}{N}
\mathbbm{E}
\left( |\mathcal{V}_N| \mathbbm{P}\left\{\left.
    {L}_{U, d_U}^{+}(N)\geq t
    \right||\mathcal{V}_N|\right\}\right) 
    dt
    \]
     \[=
\int_0^N \frac{1}{N}
\mathbbm{E}
\left( |\mathcal{V}_N| \mathbbm{E}\left\{\left.
    \mathbbm{1}\{
    {L}_{U, d_U}^{+}(N)\geq t
    \} \right||\mathcal{V}_N|\right\}\right) 
    dt.
    \]
For any $t>0$ we derive with the help of the Cauchy-Bunyakovskii inequality
\begin{equation}\notag
\mathbbm{E}\left( |\mathcal{V}_N|\mathbbm{E}\left\{\left.
    \mathbbm{1}\{
    {L}_{U, d_U}^{+}(N)\geq t
    \} \right||\mathcal{V}_N|\right\}\right)
    \end{equation} 
    \[ \leq \left(
    \mathbbm{E} |\mathcal{V}_N|^2
    \right)^{1/2}
\left( \mathbbm{E}\left( 
\mathbbm{E}\left\{\left.
    \mathbbm{1}\{
    {L}_{U, d_U}^{+}(N)\geq t
    \} \right||\mathcal{V}_N|\right\}\right)^2
    \right)^{1/2}
\]
\[ \leq \left(
    \mathbbm{E} |\mathcal{V}_N|^2
    \right)^{1/2}
\left( \mathbbm{E}
\mathbbm{E}\left\{\left.
    \mathbbm{1}\{
    {L}_{U, d_U}^{+}(N)\geq t
    \} \right||\mathcal{V}_N|\right\}
    \right)^{1/2}
\]
\[ = \left(
    \mathbbm{E} |\mathcal{V}_N|^2
    \right)^{1/2}
\left( 
\mathbbm{P}\left\{
    {L}_{U, d_U}^{+}(N)\geq t
    \right\}
    \right)^{1/2}.
\]
With the help of the uniform bound from Theorem \ref{T1} we obtain from here
\begin{equation}\notag
\mathbbm{E}\left( |\mathcal{V}_N|\mathbbm{E}\left\{\left.
    \mathbbm{1}\{
    {L}_{U, d_U}^{+}(N)\geq t
    \} \right||\mathcal{V}_N|\right\}\right)
    \end{equation} 
\[ \leq  \left(
    \mathbbm{E} |\mathcal{V}_N|^2
    \right)^{1/2} 
\left( c_2e^{-\alpha_2\sqrt{\lambda}t}
    \right)^{1/2} = \left(
    \lambda N^2 + \left(
    \lambda N^2 
    \right)^{2} 
    \right)^{1/2} 
\left( c_2e^{-\alpha_2\sqrt{\lambda}t}
    \right)^{1/2}
\]
\[\leq c_3 \lambda N^2 e^{-\frac{\alpha_2}{2} \sqrt{\lambda}t} \]
for all large $N$ and some positive constant $c_3$ independent of $\lambda$.
Substituting the last bound into (\ref{bd5}), we derive
\begin{equation}\label{bd6}    
\mathbbm{E}\left( |\mathcal{V}_N|\mathbbm{P}\left\{\left.
    {L}_{U, d_U}^{+}(N)\geq \eta_U^{+}
    \right||\mathcal{V}_N|\right\}\right)
\end{equation}
     \[\leq
\int_0^N \frac{1}{N}c_3 \lambda N^2 e^{-\frac{\alpha_2}{2} \sqrt{\lambda}t}
    dt
        \leq c_4\sqrt{\lambda}N
    \]
    for some positive constant $c_4,$ independent of $\lambda$.
By symmetry, the same bound holds for the last term in (\ref{bd2}), yielding, together with (\ref{bd6}),
\begin{equation}\label{bd7}
\mathbbm{E} \mathcal{R}_N
\leq 2c_4\sqrt{\lambda}N,
\end{equation}
which is the upper bound in Theorem \ref{th:escBound}.

\begin{figure}[h]
        \centering
        \includegraphics[width = 0.57\linewidth]{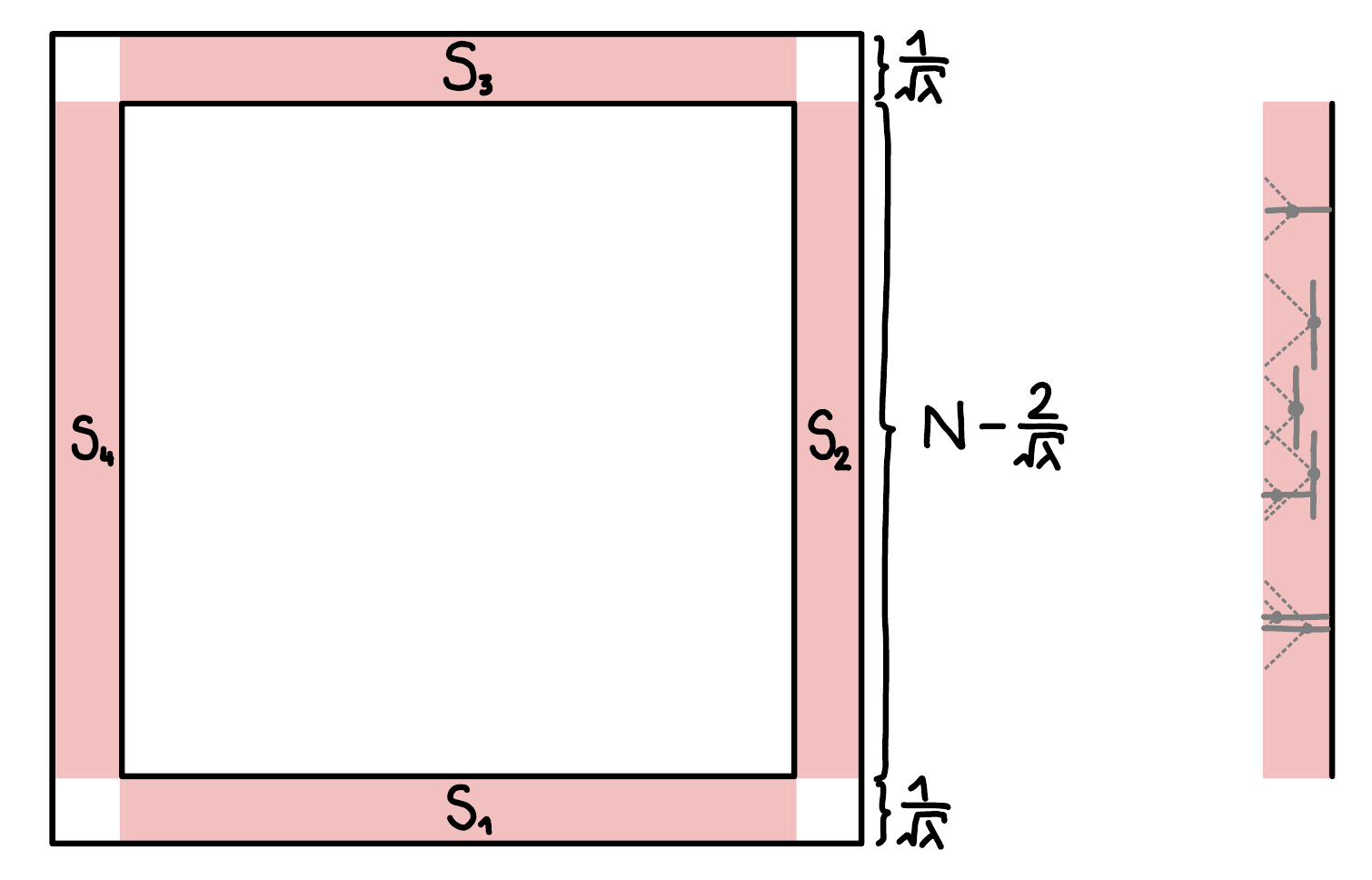}
        \caption{Illustration of the frame of width $\frac{1}{\sqrt{\lambda}}$ and stripes $S_1,\dots,S_4$. On the right side, an example of $S_2$ with tops and points that are not tops whose lines reach the border.}
        \label{fig:lowerBoundEsc}
    \end{figure}

Next, we derive the lower bound. For that, we shall simply count the escaping rays that originated only within a narrow band of width $\frac{1}{\sqrt{\lambda}}$ along the boundary of box $\Lambda_N$, as illustrated in Figure \ref{fig:lowerBoundEsc}. Cutting away the corner squares, also with a side length $\frac{1}{\sqrt{\lambda}}$, leaves the four stripes $S_1,\dots,S_4$ of equal area (shaded red in Figure \ref{fig:lowerBoundEsc}). The numbers of escaping rays 
from the four stripes are equally distributed due to symmetry. Therefore, without loss of generality, we consider stripe
$S_2$ (see Figure \ref{fig:lowerBoundEsc}) along the edge 
$$L:= \{(N,y),  \  0\leq y\leq N\}. $$
Recall Definition \ref{def:tops} of 
the notion \enquote{top}. This time we use it, 
replacing the line $\mathcal{L}_{v, d_v}^{\pm}$ with the edge $L$.
Let us rephrase it here: a point $u\in \mathcal{V}_N
$
is called a {\normalfont top} 
with respect to ${L}$
in $\mathcal{V}_N
$,
if, in order to reach the edge $L$ from this point $u$, it is necessary and sufficient to require that 
the direction vector $d_u$ is orthogonal to $L$, which is 
$d_u=(1,0)$.

By this definition, if a top within $S_2$ is assigned the direction vector  $(1,0)$, this guarantees an escaping line from this point through edge $L$. 
(This is illustrated in the example stripe $S_2$ in the Figure \ref{fig:lowerBoundEsc} as well.) 
 
    Hence, if $T$ is the number of tops in stripe $S_2$ (see  Figure \ref{fig:lowerBoundEsc})  of length $N-\frac{2}{\sqrt{\lambda}}$ and width $\frac{1}{\sqrt{\lambda}}$, then we have a rough lower bound in a form 
    \begin{equation}\label{lr}
        \mathbbm{E}\mathcal{R}_N \geq 4 \cdot\frac{1}{2}\cdot\mathbbm{E} T,
    \end{equation}
    where 4 represents 4 sides, and $\frac{1}{2}$ is the probability that a top is assigned the required direction. 
    To find now a lower bound for $\mathbbm{E} T$ in (\ref{lr}) let us
    divide the stripe $S_2$  into 
    boxes with side length $\frac{1}{\sqrt{\lambda}}$. This gives us
    \begin{equation}\label{Rn}
        n:=\left\lfloor \sqrt{\lambda}\left(N-\frac{2}{\sqrt{\lambda}}\right)\right\rfloor\geq \sqrt{\lambda}N-3
    \end{equation} 
    whole non-intersecting squares  with side length $\frac{1}{\sqrt{\lambda}}$.
   Let $M$ be the total number of tops in these $n$ boxes only, hence, $M\leq T$, and let
   $W$ be the number of non-empty boxes. As the
     number of points in each box follows a $\text{Poi}(1)$ distribution independently for all boxes, 
     the distribution of $W$ is Binomial:
    $$ W \sim\text{Bin}\left(n,1-\frac{1}{e}\right).$$ 
    As we argued in the proof of Theorem \ref{T1} (see (\ref{S5})), it holds that $$M\geq \left\lfloor \frac{W}{3}\right\rfloor-2,$$ 
    thus
    \begin{equation}\notag
     \mathbbm{E} T\geq   \mathbbm{E} M\geq \mathbbm{E}  \frac{W}{3} -3= \frac{1}{3}n\left(1-\frac{1}{e}\right)-3\geq\frac{1}{3} (\sqrt{\lambda}N-3)\left(1-\frac{1}{e}\right)-3,
    \end{equation} taking (\ref{Rn}) into account in the last inequality.
    Together with (\ref{lr}) gives us the lower bound
    \begin{equation}
        \label{lr1}
        \mathbbm{E}\mathcal{R}_N \geq 2\mathbbm{E} T  \geq \frac{2}{3}\sqrt{\lambda}\left(1-\frac{1}{e}\right)N-5.
    \end{equation}
  The lower (\ref{lr1}) and the upper (\ref{bd7}) bounds confirm the statement of the theorem. \hfill$\Box$

\section*{Acknowledgements} 
This work was supported by a fellowship from the German Academic Exchange Service (DAAD) for the author Emily Ewers.

\section*{Declarations and Statements}
The authors have no relevant financial or non-financial interests to disclose.

\end{document}